\documentclass{article}
\usepackage{amsmath,amssymb,amscd,latexsym}

\newcommand{\PreserveBackslash}[1]{\let\temp=\\#1\let\\=\temp}

\def\A{{\mathcal{A}}}      
\def\Br{{\mathfrak{B}}}    
\def\P{{\mathfrak{P}}}     
\def\RR{{\mathcal{R}}}     
\def\Gam{\Gamma}            
\def\Ar{\mathcal{A}_{\Gamma}}        
\def\pAr{\mathcal{A}_{\Gamma}^+}     
\def\pe{=_{p}}             
\def\la{{\langle}}         
\def\ra{{\rangle}}         
\def\d{{\rm deg}}          
\def\dz{{\mathcal N}}      
\def\DZ{{\mathcal N_{\Gamma}}}      

\def\P{\Bbb P}      

\renewcommand{\P}{{\mathcal P}}
\def\S{{\mathcal S}}
\def\Z{\Bbb Z}      

\def\im{\mathop{\rm im}}

\def\ker{\rm ker}
\def\im{\rm Im}

\def\Java{\hbox{\sc J\kern -1pt{}ava}}

\newtheorem{theorem}{Theorem}[section]
\newtheorem{lemma}[theorem]{Lemma}
\newtheorem{corollary}[theorem]{Corollary}
\newtheorem{proposition}[theorem]{Proposition}
\newtheorem{example}[theorem]{Example}

\newenvironment{proof}{\noindent{\bf Proof. } }{\nobreak\hfill
   \vrule height8pt width4pt\bigskip}

\renewenvironment{proof}{\noindent{\bf Proof. }}{\nobreak\hfill
    \qed
   \bigskip}

\def\Aut{{\rm Aut}}

\renewcommand{\P}{{\mathcal P}}

\def\Z{{\mathbb Z}}

\def\square{{\vcenter{\hrule height.4pt
      \hbox{\vrule width.4pt height5pt \hskip5pt
           \vrule width.4pt}
      \hrule height.4pt}}}
\def\qed{\hfill$\square$}

\begin{document}

\title{Local indicability and commutator \\ subgroups of Artin groups }
\author{Jamie Mulholland and Dale Rolfsen}
\maketitle

\begin{abstract}
 Artin groups (also known as Artin-Tits groups) are generalizations
 of Artin's braid groups. This paper concerns Artin groups of
 spherical type, that is, those whose corresponding Coxeter group
 is finite, as is the case for the braid groups. We compute
 presentations for the commutator subgroups of the irreducible
 spherical-type Artin groups, generalizing the work of Gorin and
 Lin \cite{GL69} on the braid groups.  Using these presentations we
 determine the local indicability of the irreducible spherical
 Artin groups (except for $F_4$ which at this time remains undetermined).  We end with
 a discussion of the current state of the right-orderability of
 the spherical-type Artin groups.
\end{abstract}

\section{Introduction}
\label{sec:introduction}

  A number of recent discoveries regarding the Artin braid groups
  $\Br_n$ complete a rather interesting story about the
  orderability\footnote{
  A group $G$ is
  {\it right-orderable} if there exists a
  strict total ordering $<$ of its elements which is
  right-invariant: $g < h$ implies $gk < hk$ for all $g, h, k \in
  G$.  If in addition $g < h$ implies $kg < kh$, the group is said
  to be {\it orderable}, or for emphasis, {\it
  bi-orderable}.}
  of these groups.
  These discoveries were as follows.

 In $1969$, Gorin and Lin \cite{GL69}, by computing presentations for
 the commutator subgroups $\Br_n'$ of the braid groups $\Br_n$,
 showed that $\Br_3'$ is a free group of
 rank $2$, $\Br_4'$ is the semidirect product of two free groups
 (each of rank $2$), and $\Br_n'$ is finitely generated and perfect
 for $n\ge 5$.  It follows from these results that
 $\Br_n$ is locally indicable\footnote{A group $G$ is {\it locally indicable}\index{locally indicable} if
 for every nontrivial, finitely generated subgroup $H$ of $G$ there exists a {\it nontrivial}
 homomorphism $H \to \Z$.} if and only if $n < 5$.

 Neuwirth in 1974 \cite{Ne74}, observed $\Br_n$ is {\it not} bi-orderable if $n
 \ge 3$.  Twenty years later, Dehornoy \cite{De94} showed the braid groups are in fact
 right-orderable for all $n$.  Furthermore, it has been observed, \cite{RZ98}, \cite{KRo02}, that the
 subgroups $\P_n$ of pure braids are bi-orderable.

 These orderings were fundamentally different, and it was natural to ask if there might be compatible orderings,
 that is a right-invariant ordering of $\Br_n$ which restricts to a bi-ordering of $\P_n.$
 This question was answered by Rhemtulla and Rolfsen \cite{RR02}
 by exploiting the connection between local indicability and orderability.
 They showed that since the braid groups $\Br_n$ are not locally
 indicable for $n\ge 5$ a right-ordering on $\Br_n$ could not restrict to a
 bi-ordering on $\P_n$ (or on any subgroup of finite index).

 This paper is concerned with investigating the extent to which of these results on the
 braid groups extend to other Artin groups, or at least those of spherical type (defined in the next section).
 In  particular, we are concerned with determining the local indicability of the spherical Artin groups.
 Because the full details of the Gorin-Lin calculations do not seem to appear in the literature, we present
 a fairly comprehensive account of the calculation of commutator subgroups of the braid groups, which are the
 Artin groups of type $A_n$.  These methods, essentially the Reidemeister-Schreier method plus a few tricks, are
 also used to calculate presentations of the commutator subgroups of the other spherical Artin groups.

 In the next section we will define Coxeter graphs $\Gamma$, the corresponding Coxeter groups $W_\Gamma$ and the Artin
 groups $\A_\Gamma$.   The spherical Artin groups are classified according to types:
 $A_n (n \ge 1), B_n (n \ge 2), D_n (n \ge 4), E_6, E_7, E_8, F_4, H_3, H_4$ and $I_2(m), (m \ge 5)$.
 Our main results are summarized in the following theorem, where $\A_\Gamma '$ denotes the commutator subgroup.
 Recall that a perfect group is one which equals its own commutator subgroup; any homomorphism from a perfect
 group to an abelian group must be trivial.

\begin{theorem}
  The following commutator subgroups are finitely generated and perfect:
   \begin{itemize}
    \item[1.] $\A_{A_n}'$ for $n\ge 4$,
    \item[2.] $\A_{B_n}'$ for $n\ge 5$,
    \item[3.] $\A_{D_n}'$ for $n\ge 5$,
    \item[4.] $\A_{E_n}'$ for $n=6,7,8$,
    \item[5.] $\A_{H_n}'$ for $n=3,4$.
   \end{itemize}
Hence, the corresponding Artin groups are not locally indicable.
\end{theorem}

On the other hand, we show the remaining spherical-type Artin groups
{\it are} locally indicable (excluding the type $F_4$ which at
this time remains undetermined) .

In a final section we discuss the orderability of
the spherical-type Artin groups.   We show that to determine
the orderability of the spherical-type Artin
groups it is sufficient to consider the positive Artin
monoid. Furthermore, we show
that to prove {\it all} spherical-type Artin groups are
right-orderable it would suffice to show the Artin group (or monoid) of type $E_8$
is right-orderable.

 \section{Coxeter and Artin groups}
 \label{sec:coxeter-definition}

 Let $S$ be a finite set.  A {\bf Coxeter matrix} \index{Coxeter!matrix} over $S$ is a
 matrix $M=(m_{ss'})_{s,s'\in S}$ indexed by the elements of $S$
 and satisfying
 \begin{flushleft}
 ~~~~~(a)~~$m_{ss}=1$ if $s \in S$, \\
 ~~~~~(b)~~$m_{ss'}=m_{s's} \in \{2, \ldots , \infty \}$ if
 $s,s' \in S$ and $s \neq s'$.
 \end{flushleft}
 A Coxeter matrix $M=(m_{ss'})_{s,s'\in S}$ is usually
 represented by its {\bf Coxeter graph} \index{Coxeter!graph} $\Gamma$.  This is defined
 by the following data.
 \begin{flushleft}
   ~~~~~(a)~~$S$ is the set of vertices of $\Gamma$. \\
   ~~~~~(b)~~Two vertices $s,s' \in S$ are joined by an edge if $m_{ss'} \geq 3$. \\
   ~~~~~(c)~~The edge joining two vertices $s,s' \in S$ is labelled
   by $m_{ss'}$ if $m_{ss'}\geq 4$.
 \end{flushleft}
 The {\bf Coxeter system} \index{Coxeter!system} of type $\Gamma$ (or $M$)
 is the pair ($W,S$) where $W$ is the group having the
 presentation
 \begin{eqnarray*}
   W = \langle s \in S : (ss')^{m_{ss'}}=1 \text{ if } m_{ss'} <
   \infty \rangle.
 \end{eqnarray*}

\begin{example}
\label{ex:symmetric-group}
 \index{group!symmetric}
 It is well known that the symmetric group on $(n+1)$-letters is the
 Coxeter group associated with the Coxeter graph;\\
  \begin{picture}(300,25)(-50,-10)
    \linethickness{2pt}
     \put(0,5){\circle*{6}}
     \put(30,5){\circle*{6}}
     \put(60,5){\circle*{6}}
     \put(90,5){\circle*{3}}
     \put(100,5){\circle*{3}}
     \put(110,5){\circle*{3}}
     \put(140,5){\circle*{6}}
     \put(170,5){\circle*{6}}
     \put(200,5){\circle*{6}}
    \put(0,5){\line(1,0){60}}
    \put(140,5){\line(1,0){60}}
    \put(0,0){\makebox(0,0)[t]{\scriptsize$1$}}
    \put(30,0){\makebox(0,0)[t]{\scriptsize$2$}}
    \put(60,0){\makebox(0,0)[t]{\scriptsize$3$}}
    \put(140,0){\makebox(0,0)[t]{\scriptsize$n-2$}}
    \put(170,0){\makebox(0,0)[t]{\scriptsize$n-1$}}
    \put(200,0){\makebox(0,0)[t]{\scriptsize$n$}}
   \end{picture} \\
 \noindent where vertex $i$ corresponds to the transposition
$(i~i+1)$.
\end{example}

  If ($W,S$) is a Coxeter system with Coxeter graph $\Gamma$ (resp.
 Coxeter matrix $M$) then we say that $\Gamma$ (resp. $M$) is of
 {\bf spherical-type} if $W$ is finite.  If $\Gamma$ is connected, then $W$ is said to be
 {\bf irreducible}.  Coxeter  \cite{Co34} classified all irreducible Coxeter groups which are
 finite, a result that plays a central role in the theory of Lie groups.  We refer the reader to
 \cite{Hu72} (see also Bourbaki \cite{Bo72}, \cite{Bo02}) for further details on Coxeter groups, including a
 proof of the following.

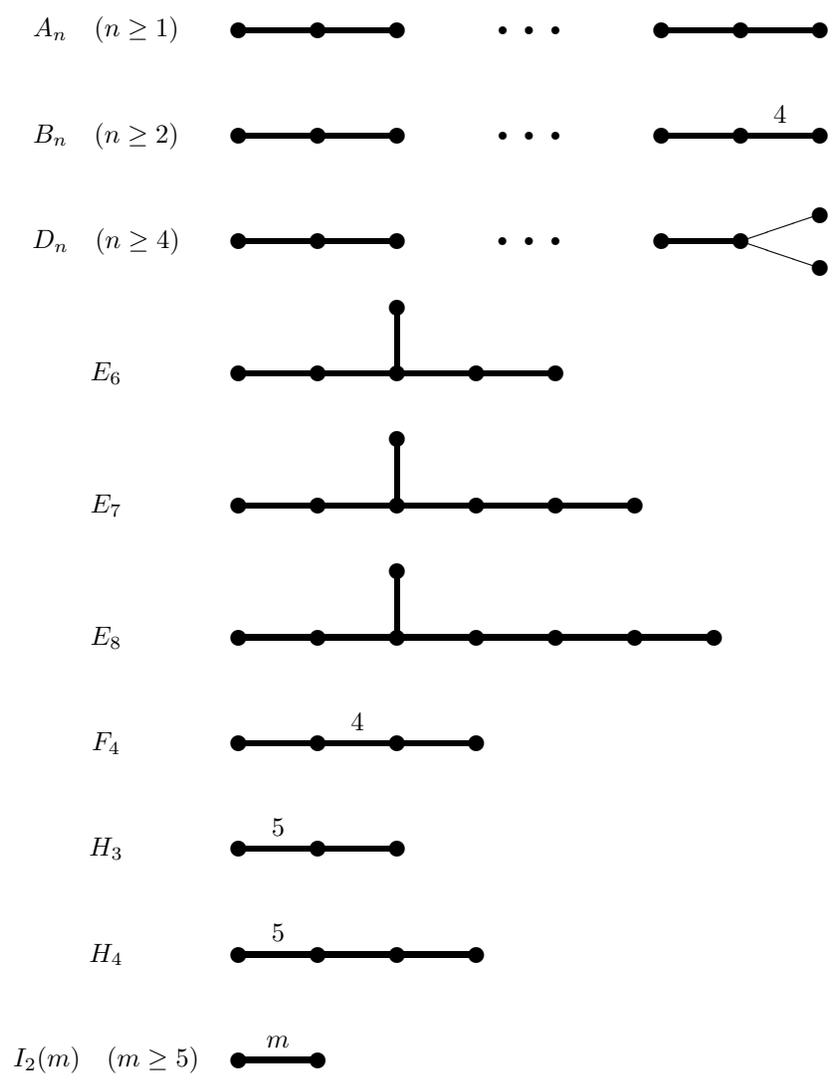
\begin{figure}[hp]
%

 \begin{picture}(300,400)(-100,-10)
  \linethickness{2pt}
    \put(-50,400){\makebox(0,0)[c]{$A_n~~~(n \geq 1)$}}
    \put(0,400){\circle*{6}}
    \put(30,400){\circle*{6}}
    \put(60,400){\circle*{6}}
    \put(100,400){\circle*{3}}
    \put(110,400){\circle*{3}}
    \put(120,400){\circle*{3}}
    \put(160,400){\circle*{6}}
    \put(190,400){\circle*{6}}
    \put(220,400){\circle*{6}}
    \put(0,400){\line(1,0){60}}
    \put(160,400){\line(1,0){60}}

    \put(-50,360){\makebox(0,0)[c]{$B_n~~~(n \geq 2)$}}
    \put(0,360){\circle*{6}}
    \put(30,360){\circle*{6}}
    \put(60,360){\circle*{6}}
    \put(100,360){\circle*{3}}
    \put(110,360){\circle*{3}}
    \put(120,360){\circle*{3}}
    \put(160,360){\circle*{6}}
    \put(190,360){\circle*{6}}
    \put(220,360){\circle*{6}}
    \put(0,360){\line(1,0){60}}
    \put(160,360){\line(1,0){60}}
    \put(205,365){\makebox(0,0)[b]{$4$}}

    \put(-50,320){\makebox(0,0)[c]{$D_n~~~(n\geq 4)$}}
    \put(0,320){\circle*{6}}
    \put(30,320){\circle*{6}}
    \put(60,320){\circle*{6}}
    \put(100,320){\circle*{3}}
    \put(110,320){\circle*{3}}
    \put(120,320){\circle*{3}}
    \put(160,320){\circle*{6}}
    \put(190,320){\circle*{6}}
    \put(220,330){\circle*{6}}
    \put(220,310){\circle*{6}}
    \put(0,320){\line(1,0){60}}
    \put(160,320){\line(1,0){30}}
    \put(190,320){\line(3,1){30}}
    \put(190,320){\line(3,-1){30}}

    \put(-50,270){\makebox(0,0)[c]{$E_6$}}
    \put(0,270){\circle*{6}}
    \put(30,270){\circle*{6}}
    \put(60,270){\circle*{6}}
    \put(90,270){\circle*{6}}
    \put(120,270){\circle*{6}}
    \put(60,295){\circle*{6}}
    \put(0,270){\line(1,0){120}}
    \put(60,270){\line(0,1){25}}

    \put(-50,220){\makebox(0,0)[c]{$E_7$}}
    \put(0,220){\circle*{6}}
    \put(30,220){\circle*{6}}
    \put(60,220){\circle*{6}}
    \put(90,220){\circle*{6}}
    \put(120,220){\circle*{6}}
    \put(150,220){\circle*{6}}
    \put(60,245){\circle*{6}}
    \put(0,220){\line(1,0){150}}
    \put(60,220){\line(0,1){25}}

    \put(-50,170){\makebox(0,0)[c]{$E_8$}}
    \put(0,170){\circle*{6}}
    \put(30,170){\circle*{6}}
    \put(60,170){\circle*{6}}
    \put(90,170){\circle*{6}}
    \put(120,170){\circle*{6}}
    \put(150,170){\circle*{6}}
    \put(180,170){\circle*{6}}
    \put(60,195){\circle*{6}}
    \put(0,170){\line(1,0){180}}
    \put(60,170){\line(0,1){25}}

    \put(-50,130){\makebox(0,0)[c]{$F_4$}}
    \put(0,130){\circle*{6}}
    \put(30,130){\circle*{6}}
    \put(60,130){\circle*{6}}
    \put(90,130){\circle*{6}}
    \put(0,130){\line(1,0){90}}
    \put(45,135){\makebox(0,0)[b]{$4$}}

    \put(-50,90){\makebox(0,0)[c]{$H_3$}}
    \put(0,90){\circle*{6}}
    \put(30,90){\circle*{6}}
    \put(60,90){\circle*{6}}
    \put(0,90){\line(1,0){60}}
    \put(15,95){\makebox(0,0)[b]{$5$}}

    \put(-50,50){\makebox(0,0)[c]{$H_4$}}
    \put(0,50){\circle*{6}}
    \put(30,50){\circle*{6}}
    \put(60,50){\circle*{6}}
    \put(90,50){\circle*{6}}
    \put(0,50){\line(1,0){90}}
    \put(15,55){\makebox(0,0)[b]{$5$}}

    \put(-50,10){\makebox(0,0)[c]{$I_2(m)~~~(m \geq 5)$}}
    \put(0,10){\circle*{6}}
    \put(30,10){\circle*{6}}
    \put(0,10){\line(1,0){30}}
    \put(15,15){\makebox(0,0)[b]{$m$}}

\end{picture}
  \caption{The connected Coxeter graphs of spherical type}
  \label{fig:coxeter-graphs}
 \end{figure}

 \begin{theorem}  The connected Coxeter graphs of spherical type are exactly those depicted in figure~\ref{fig:coxeter-graphs}.
 \end{theorem}

 The letter beside each of the graphs in figure~\ref{fig:coxeter-graphs}
 is called the {\bf type}  of
 the Coxeter graph; the subscript denotes the number of
 vertices.  Recall example~\ref{ex:symmetric-group} shows
 the symmetric group on (n+1)-letters is a Coxeter group of type
 $A_n$.

Let $M$ be a Coxeter matrix over $S$ as described above, and let $\Gamma$ be the
corresponding Coxeter graph.  Fix a set $\Sigma$ in one-to-one
correspondence with $S$.

We adopt the following notation, where $a, b \in \Sigma$:
\begin{eqnarray*}
  \langle ab \rangle^q=\underbrace{aba \ldots}_{\text{q factors}}
\end{eqnarray*}

 The {\bf Artin system} of type $\Gamma$ (or $M$) is the pair
($\A,\Sigma$) where $\A$ is the group having presentation
\begin{eqnarray*}
  \A = \langle a \in \Sigma : \langle ab \rangle^{m_{ab}}=\langle ba
  \rangle^{m_{ab}} \text{ if } m_{ab}< \infty \rangle.
\end{eqnarray*}
The group $\A$ is called the {\bf Artin group} \index{Artin!group}
of type $\Gamma$ (or $M$), and is sometimes denoted by $\Ar$. So,
as with Coxeter systems, an Artin system is an Artin group with
a prescribed set of generators.

There is a natural map $\nu : \Ar \longrightarrow W_{\Gamma}$
sending generator $a_i \in \Sigma$ to the corresponding generator
$s_i \in \S$.  This map is indeed a homomorphism since the
equation $\langle s_is_j \rangle^{m_{ij}}=\langle s_js_i
\rangle^{m_{ij}}$ follows from $s_i^2=1$, $s_j^2=1$ and
$(s_is_j)^{m_{ij}}=1$.  Since $\nu$ is clearly surjective it
follows that the Coxeter group $W_{\Gamma}$ is a quotient of the
Artin group $\Ar$.  The kernel of $\nu$ is called the {\bf pure
Artin group}, generalizing the
definition of the pure braid group.

\begin{example}  The Artin group $\A_{A_n}$ is isomorphic with the braid group
$\Br_{n+1}$ on $n+1$ strings.  The homomorphism $\A_{A_n} \to W_{A_n}$ corresponds to the map which assigns to each braid the permutation determined by running from one end of the strings to the other.  Pure braids are those for which the permutation is the identity.
\end{example}

 The Artin group of a spherical-type Coxeter graph is called an Artin
 group of {\bf spherical-type}, that is, the corresponding
 Coxeter group $W_{\Gam}$ is finite.  An Artin group $\Ar$ is
 called {\bf irreducible} if the Coxeter
 graph $\Gam$ is connected. In particular, the Artin groups
 corresponding to the graphs in figure~\ref{fig:coxeter-graphs} are those which are
 irreducible and of spherical-type. These Artin groups are by far the most well-understood,
 and are our main interest in the remaining sections.

 Van der Lek \cite{Le83} has shown that for each subgraph $I \subset \Sigma$ the
 corresponding subgroup and subgraph are an Artin system.
 That is, parabolic subgroups of Artin groups (those generated by a subset of the generators) are indeed Artin groups.
 A proof of this fact also appears in \cite{Pa97}.  Thus inclusions
 among Coxeter graphs give rise to inclusions for the associated Artin
 groups.  Crisp \cite{Cr99} shows quite a
 few more inclusions hold among the irreducible spherical-type Artin
 groups.  Table~\ref{tab:artin-injections} summarizes these inclusions.
 Notice that every irreducible spherical-type Artin group embeds into
 an Artin group of type $A$, $D$ or $E$.

\begin{table}[top]
  \begin{center}
    \begin{tabular}{|c|c|} \hline
      \multicolumn{2}{|c|}{$\Ar$ injects into $A_{\Gamma'}$} \\ \hline \hline
         $\Gamma$    & $\Gamma'$                   \\ \hline \hline
         $A_n$       & $A_{m}$ ($m\ge n$), \\
                     & $B_{n+1}$ ($n\ge 2$), \\
                     & $D_{n+2}$, \\
                     & $E_6$ ($ 1\le n\le 5$), \\
                     & $E_7$ ($ 1\le n\le 6$), \\
                     & $E_8$ ($ 1\le n\le 7$), \\
                     & $F_4$, $H_3$  ($ 1\le n\le 2$), \\
                     & $H_4$ ($ 1\le n\le 3$)  \\
                     & $I_{2}(3)$ ($ 1\le n\le 2$)    \\  \hline
         $B_n$       & $A_{n}$, $A_{2n-1}$, $A_{2n}$, $D_{n+1}$    \\ \hline
         $E_6$       & $E_7$, $E_8$         \\ \hline
         $E_7$       & $E_8$         \\ \hline
         $F_4$       & $E_6$, $E_7$, $E_8$ \\ \hline
         $H_3$       & $D_6$         \\ \hline
         $H_4$       & $E_8$         \\ \hline
         $I_{2}(m)$  & $A_{m-1}$     \\ \hline
    \end{tabular}
  \end{center}
  \caption{Inclusions among Artin groups}
  \label{tab:artin-injections}
 \end{table}

 Cohen and Wales \cite{CW01} use the fact that irreducible finite
 type Artin groups embed into an Artin group of type $A$, $D$ or $E$
 to show all Artin groups of spherical-type are {\bf linear}
 (have a faithful finite-dimensional linear representation) by showing Artin groups of
 type $D$, and $E$ are linear, thus generalizing the recent result
 that the braid groups (Artin groups of type $A$) are linear
 \cite{Bi01}, \cite{Kra02}.

 We close this section by noting that Deligne \cite{De72} showed that each Artin group of
 spherical-type appears as the fundamental group of the complement of a complex hyperplane arrangement,
 which is an Eilenberg-Maclane space.  From this
 point of view we can see that spherical-type Artin groups are torsion free and have finite cohomological dimension.

\section{Commutator Subgroups}

 Our basic tool for finding presentations for the commutator subgroups of Artin groups
 is the classical Reidemeister-Schreier method.  To fix notation, we begin with a brief
 review of this algorithm.   For a more complete discussion, see \cite{MKS76}.

\subsection{Reidemeister-Schreier algorithm}
\label{sec:reidemeister-schreier-method}

  Let $G$ be an arbitrary group with presentation
  $\la a_{1}, \ldots,a_{n} : R_{\mu}(a_{\nu}), \ldots \ra$ and $H$ a subgroup of $G$.
  A system of words $\RR$ in the generators $a_{1}, \ldots, a_{n}$ is called a
  {\bf Schreier system}\index{Schreier-system} for G modulo H
  if (i) every right coset of $H$ in $G$
  contains exactly one word of $\RR$ (i.e. $\RR$ forms
  a system of right coset representatives), (ii) for each word in $\RR$ any
  initial segment is also in $\RR$ (i.e. initial segments of right coset
  representatives are again right coset representatives).  Such a Schreier system
  always exists, see for example \cite{MKS76}.  Suppose now that we have fixed a Schreier
  system $\RR$.  For each word $W$ in the generators $a_1, \ldots, a_n$ we let
  $\overline{W}$ denote the unique representative in $\RR$
  of the right coset $HW$.  Denote

  \begin{eqnarray}
     s_{K,a_v} = Ka_v \cdot \overline{Ka_{v}}^{-1},
  \end{eqnarray}
  \noindent
  for each $K \in \RR$ and generator $a_v$.  A theorem of Reidemeister-Schreier
  (theorem $2.9$ in [MKS76]) states that
  $H$ has presentation
  \begin{eqnarray}{\label{general-presentation}}
     \la s_{K,a_{\nu}}, \ldots : s_{M,a_{\lambda}} , \ldots, \tau(KR_{\mu}K^{-1}),\ldots
     \ra
  \end{eqnarray}
  \noindent
  where $K$ is an arbitrary  Schreier representative, $a_v$ is an
  arbitrary generator and $R_{\mu}$ is an arbitrary defining
  relator in the presentation of $G$, and $M$ is a Schreier
  representative and $a_{\lambda}$ a generator such that
  \begin{eqnarray*}
     Ma_{\lambda} \thickapprox \overline{Ma_{\lambda}},
  \end{eqnarray*}
  \noindent
  where $\thickapprox$ means "freely equal"\index{freely equal}, i.e. equal in the free group generated
  by \linebreak $\{a_1, \ldots, a_n \}$.  The function $\tau$ is
  a {\bf Reidemeister rewriting function}\index{Reidemeister!rewriting function} and is defined
  according to the rule
  \begin{eqnarray}{\label{reidemeister-rewriting-function}}
    \tau(a_{i_1}^{\epsilon_1}\cdots a_{i_p}^{\epsilon_p})=s_{K_{i_1},a_{i_1}}^{\epsilon_1}\cdots s_{K_{i_p},a_{i_p}}^{\epsilon_p}
  \end{eqnarray}
  \noindent
  where $K_{i_j}=\overline{a_{i_1}^{\epsilon_1}\cdots
  a_{i_{j-1}}^{\epsilon_{j-1}}}$, if $\epsilon_j = 1$, and $K_{i_j}=\overline{a_{i_1}^{\epsilon_1}\cdots
  a_{i_{j}}^{\epsilon_{j}}}$, if $\epsilon_j=-1$.  It should be
  noted that computation of $\tau(U)$ can be carried out by
  replacing a symbol $a_{v}^{\epsilon}$ of U by the appropriate
  s-symbol $s_{K,a_{\nu}}^{\epsilon}$.  The main property of a
  Reidemeister rewriting function is that for an element $U \in H$
  given in terms of the generators $a_{\nu}$ the word $\tau(U)$ is
  the same element of $H$ rewritten in terms of the generators
  $s_{K,a_{\nu}}$.

\subsection{Characterization of the commutator subgroups}
\label{sec:characterization}

The {\bf commutator subgroup}\index{commutator subgroup} $G'$ of a
group $G$ is the subgroup generated by the elements $[g_1,g_2] :=
g_{1}g_{2}g_{1}^{-1}g_{2}^{-1}$ for all $g_1,g_2 \in G$.  Such
elements are called {\bf commutators}\index{commutator}. It is an
elementary fact in group theory that $G'$ is a normal subgroup in
$G$ and the quotient group $G/G'$ is abelian.  In fact, for any
normal subgroup $N \lhd G$ the quotient group $G/N$ is abelian if
and only if $G' < N$.  If $G$ is given in terms of a presentation
$\la \mathcal{G} : \RR \ra$ where $\mathcal {G}$ is a set of
generators and $\RR$ is a set of relations, then a presentation
for $G/G'$ is obtained by abelianizing the
presentation for $G$, that is, by adding relations $gh=hg$ for all
$g,h \in \mathcal{G}$. This is denoted by $\la \mathcal{G} : \RR
\ra_{\mathrm{Ab}}$.

Let $U \in \Ar$, and write $U=a_{i_1}^{\epsilon_1} \cdots
a_{i_r}^{\epsilon_r}$, where $\epsilon_i = \pm 1$. The ({\bf
canonical}) {\bf degree of}\index{Artin!group!degree} $U$ is
defined to be
\begin{center}
  $\d(U) := \sum_{j=1}^{r}\epsilon_{j}$.
\end{center}
Since each defining relator in the presentation for $\Ar$ has
degree equal to zero the map $\d$ is a well defined homomorphism
from $\Ar$ into $\Z$. Let $\DZ$ denote the kernel of $\d$; $\DZ =
\{ U \in \Ar : \d(U) =0 \}$.  It is a well known fact that for the
braid group (i.e. $\Gamma = A_n$) $\dz_{A_n}$ is precisely the
commutator subgroup. In this section we generalize this fact for
all Artin groups.

Let $\Gam_{odd}$ denote the graph obtained from $\Gam$ by removing
all the even-labelled edges and the edges labelled $\infty$. The
following theorem tells us exactly when the commutator subgroup
$\Ar'$ is equal to $\DZ$.

\begin{proposition}
\label{thm:connected-commutator}
  For an Artin group $\Ar$, $\Gam_{odd}$ is
  connected if and only if the commutator subgroup $\Ar'$ is equal
  to $\DZ$.
\end{proposition}

\begin{proof}
 For the direction $(\Longrightarrow)$
 the connectedness of $\Gamma_{odd}$ implies
 \begin{center}
   $\Ar/\Ar' \simeq \Z$.
 \end{center}
 Indeed, start with any generator $a_i$, for any other generator
 $a_j$ there is a path from $a_i$ to $a_j$ in $\Gam_{odd}$:
 \begin{center}
   $a_i=a_{i_i} \longrightarrow a_{i_2} \longrightarrow \cdots
   \longrightarrow a_{i_m} = a_j$.
 \end{center}
 Since $m_{i_{k}i_{k+1}}$ is odd the relation
 \begin{center}
    $\la a_{i_k}a_{i_{k+1}} \ra^{m_{i_{k}i_{k+1}}} = \la a_{i_{k+1}}a_{i_{k}} \ra^{m_{i_{k}i_{k+1}}}$
 \end{center}
 becomes $a_{i_k}=a_{i_{k+1}}$ in
 $\Ar/\Ar'$.  Hence, $a_i=a_j$ in $\Ar/\Ar'$.  It follows that,
 \begin{eqnarray*}
   \Ar/\Ar' &\simeq& \la a_1,\ldots,a_n : a_1=\cdots=a_n \ra \\
             &\simeq& \Z,
 \end{eqnarray*}
 where the isomorphism $\phi : \Ar/\Ar' \longrightarrow \Z$ is given by
 \begin{center}
   $U\Ar' \longmapsto \d(U)$.
 \end{center}
 Therefore, $\Ar' = \ker \phi = \DZ$.

 We leave the proof of the other direction to proposition
 \ref{thm:characterization-commuator-subgroup}, where a more
 general result is stated.
\end{proof}

For the case when $\Gam_{odd}$ is not connected we can get a more
general description of $\Ar'$ as follows.  Let $\Gam_{odd}$ have
$m$ connected components; $\Gam_{odd} = \Gam_1 \sqcup \ldots
\sqcup \Gam_m$.  Let $\Sigma_i \subset \Sigma$ be the
corresponding sets of vertices. For each $1\le k \le m$ define the
map
 \begin{center}
   $\d_k : \Ar \longrightarrow \Z$
 \end{center}
as follows:  If $U=a_{i_1}^{\epsilon_1} \cdots
a_{i_r}^{\epsilon_r} \in \Ar$ take
 \begin{equation*}
   \d_k(U) = \sum_{1\le j \le r \; \mathbf{where}\ a_{i_j}\in \Sigma_k}
   \epsilon_j.
 \end{equation*}
 It is straight forward to check that
 for each $1\le k \le m$ the map $\d_k$ agrees on $\la ab
 \ra^{m_{ab}}$ and $\la ba \ra^{m_{ab}}$ for all $a,b \in \Sigma$.
 Hence, $\d_k : \Ar \longrightarrow \Z$ is a homomorphism for each
 $1\le k \le m$.  We combine these $m$ degree maps to get the
 following homomorphism:
 $$ \d_{\Gam} : \Ar \longrightarrow \Z^m $$ by
 $$\d_{\Gam}(U)=(\d_1(U), \ldots, \d_m(U)).$$
 When $\Gam_{odd}$ is connected, i.e. $m=1$, $\d_{\Gam}$ is just
 the canonical degree.  For $U\in \Ar$ we call $\d_{\Gam}(U)$ the
 {\bf degree} of $U$.
 The following theorem tells us that the kernel of $\d_{\Gam}$ is
 precisely the commutator subgroup of $\Ar$.

\begin{proposition}
 \label{thm:characterization-commuator-subgroup}
   Let $\Gam$ be a Coxeter graph such that $\Gam_{odd}$ has $m$
   connected components.  Then $\Ar' = \ker(\d_{\Gam})$ and $\Ar / \Ar' \cong \Z^m.$
 \end{proposition}

 \begin{proof}
   This follows from
  \begin{eqnarray*}
     \Ar/ \Ar' &\simeq& \la a_1, \ldots, a_n :
              \la a_ia_j\ra^{m_{a_ia_j}} = \la a_ja_i \ra^{m_{a_ia_j}} \ra_{\mathrm{Ab}} \\
     &\simeq& \la a_1, \ldots, a_n :  a_i =a_j \text{~~iff $i$
     and j lie in the same connected} \\
     &~& ~~~~~~~~~~~~~~~~~~~~~~\text{component of $\Gam_{odd}$} \ra_{\mathrm{Ab}}, \\
     &\simeq& \Z^m,
  \end{eqnarray*}
  with the isomorphism given by
  $$U\Ar' \longmapsto (\d_1(U), \ldots,\d_m(U))=\d_{\Gam}(U),$$
  In other words, $\d_{\Gam}$ is precisely the 
  abelianization map on $\Ar$.
  \end{proof}

  It follows that $\Ar$ and $\Ar'$ fit into a short exact sequence
  $$
 \begin{CD}
  1 \longrightarrow \Ar' \longrightarrow \Ar @>{\d_{\Gam}}>> \Z^m
  \longrightarrow 1,
 \end{CD}
 $$

%
%

\subsection{Computing the presentations}
\label{sec:computing-presentations}

In this section we compute presentations for the commutator
subgroups of the irreducible spherical-type Artin groups.  We will
show that, for the most part, the commutator subgroups are
finitely generated and perfect  (equal
to its commutator subgroup).

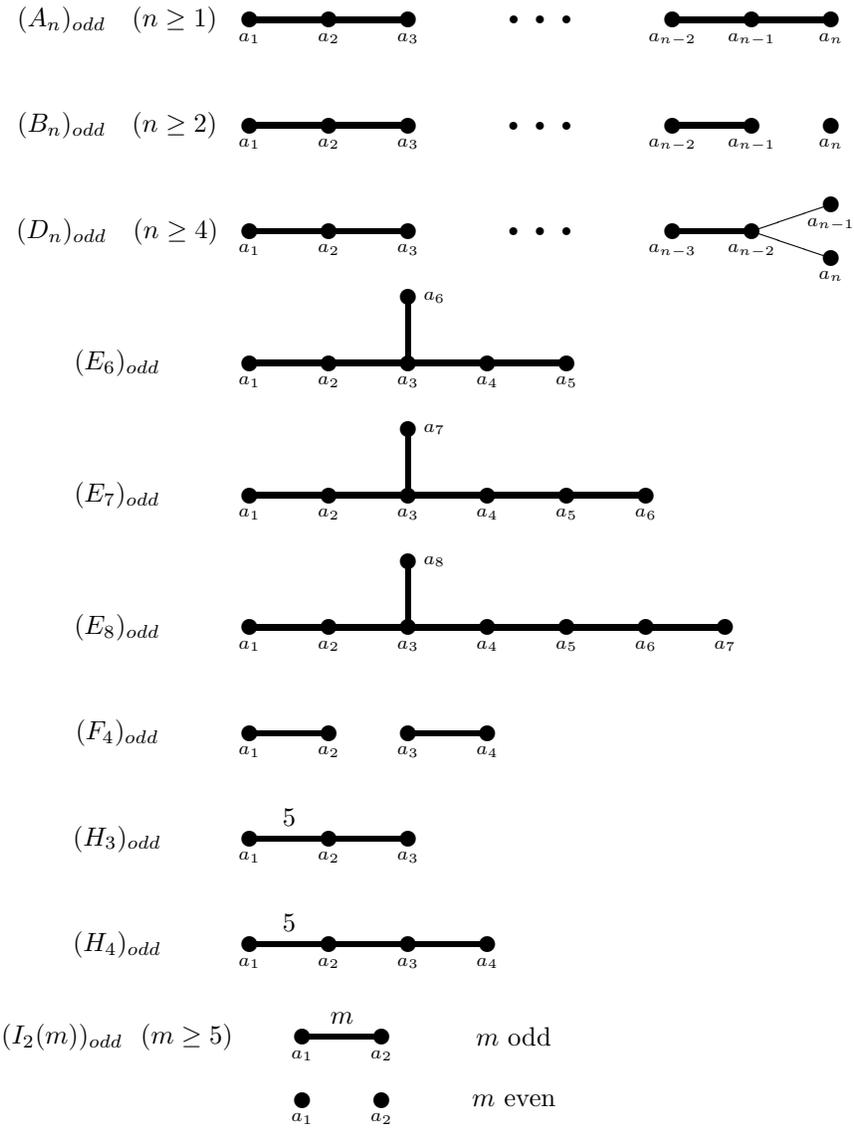
\begin{figure}[hp]
%

 \begin{picture}(300,400)(-100,-15)
  \linethickness{2pt}
    \put(-50,400){\makebox(0,0)[c]{$(A_n)_{odd}~~~(n \geq 1)$}}
    \put(0,400){\circle*{6}}
    \put(30,400){\circle*{6}}
    \put(60,400){\circle*{6}}
    \put(100,400){\circle*{3}}
    \put(110,400){\circle*{3}}
    \put(120,400){\circle*{3}}
    \put(160,400){\circle*{6}}
    \put(190,400){\circle*{6}}
    \put(220,400){\circle*{6}}
    \put(0,400){\line(1,0){60}}
    \put(160,400){\line(1,0){60}}
    \put(0,395){\makebox(0,0)[t]{\scriptsize$a_1$}}
    \put(30,395){\makebox(0,0)[t]{\scriptsize$a_2$}}
    \put(60,395){\makebox(0,0)[t]{\scriptsize$a_3$}}
    \put(160,395){\makebox(0,0)[t]{\scriptsize$a_{n-2}$}}
    \put(190,395){\makebox(0,0)[t]{\scriptsize$a_{n-1}$}}
    \put(220,395){\makebox(0,0)[t]{\scriptsize$a_n$}}

    \put(-50,360){\makebox(0,0)[c]{$(B_n)_{odd}~~~(n \geq 2)$}}
    \put(0,360){\circle*{6}}
    \put(30,360){\circle*{6}}
    \put(60,360){\circle*{6}}
    \put(100,360){\circle*{3}}
    \put(110,360){\circle*{3}}
    \put(120,360){\circle*{3}}
    \put(160,360){\circle*{6}}
    \put(190,360){\circle*{6}}
    \put(220,360){\circle*{6}}
    \put(0,360){\line(1,0){60}}
    \put(160,360){\line(1,0){30}}
    \put(0,355){\makebox(0,0)[t]{\scriptsize$a_1$}}
    \put(30,355){\makebox(0,0)[t]{\scriptsize$a_2$}}
    \put(60,355){\makebox(0,0)[t]{\scriptsize$a_3$}}
    \put(160,355){\makebox(0,0)[t]{\scriptsize$a_{n-2}$}}
    \put(190,355){\makebox(0,0)[t]{\scriptsize$a_{n-1}$}}
    \put(220,355){\makebox(0,0)[t]{\scriptsize$a_n$}}

    \put(-50,320){\makebox(0,0)[c]{$(D_n)_{odd}~~~(n\geq 4)$}}
    \put(0,320){\circle*{6}}
    \put(30,320){\circle*{6}}
    \put(60,320){\circle*{6}}
    \put(100,320){\circle*{3}}
    \put(110,320){\circle*{3}}
    \put(120,320){\circle*{3}}
    \put(160,320){\circle*{6}}
    \put(190,320){\circle*{6}}
    \put(220,330){\circle*{6}}
    \put(220,310){\circle*{6}}
    \put(0,320){\line(1,0){60}}
    \put(160,320){\line(1,0){30}}
    \put(190,320){\line(3,1){30}}
    \put(190,320){\line(3,-1){30}}
    \put(0,315){\makebox(0,0)[t]{\scriptsize$a_1$}}
    \put(30,315){\makebox(0,0)[t]{\scriptsize$a_2$}}
    \put(60,315){\makebox(0,0)[t]{\scriptsize$a_3$}}
    \put(160,315){\makebox(0,0)[t]{\scriptsize$a_{n-3}$}}
    \put(190,315){\makebox(0,0)[t]{\scriptsize$a_{n-2}$}}
    \put(220,325){\makebox(0,0)[t]{\scriptsize$a_{n-1}$}}
    \put(220,305){\makebox(0,0)[t]{\scriptsize$a_n$}}

    \put(-50,270){\makebox(0,0)[c]{$(E_6)_{odd}$}}
    \put(0,270){\circle*{6}}
    \put(30,270){\circle*{6}}
    \put(60,270){\circle*{6}}
    \put(90,270){\circle*{6}}
    \put(120,270){\circle*{6}}
    \put(60,295){\circle*{6}}
    \put(0,270){\line(1,0){120}}
    \put(60,270){\line(0,1){25}}
    \put(0,265){\makebox(0,0)[t]{\scriptsize$a_1$}}
    \put(30,265){\makebox(0,0)[t]{\scriptsize$a_2$}}
    \put(60,265){\makebox(0,0)[t]{\scriptsize$a_3$}}
    \put(90,265){\makebox(0,0)[t]{\scriptsize$a_4$}}
    \put(120,265){\makebox(0,0)[t]{\scriptsize$a_5$}}
    \put(70,295){\makebox(0,0)[c]{\scriptsize$a_6$}}

    \put(-50,220){\makebox(0,0)[c]{$(E_7)_{odd}$}}
    \put(0,220){\circle*{6}}
    \put(30,220){\circle*{6}}
    \put(60,220){\circle*{6}}
    \put(90,220){\circle*{6}}
    \put(120,220){\circle*{6}}
    \put(150,220){\circle*{6}}
    \put(60,245){\circle*{6}}
    \put(0,220){\line(1,0){150}}
    \put(60,220){\line(0,1){25}}
    \put(0,215){\makebox(0,0)[t]{\scriptsize$a_1$}}
    \put(30,215){\makebox(0,0)[t]{\scriptsize$a_2$}}
    \put(60,215){\makebox(0,0)[t]{\scriptsize$a_3$}}
    \put(90,215){\makebox(0,0)[t]{\scriptsize$a_4$}}
    \put(120,215){\makebox(0,0)[t]{\scriptsize$a_5$}}
    \put(150,215){\makebox(0,0)[t]{\scriptsize$a_6$}}
    \put(70,245){\makebox(0,0)[c]{\scriptsize$a_7$}}

    \put(-50,170){\makebox(0,0)[c]{$(E_8)_{odd}$}}
    \put(0,170){\circle*{6}}
    \put(30,170){\circle*{6}}
    \put(60,170){\circle*{6}}
    \put(90,170){\circle*{6}}
    \put(120,170){\circle*{6}}
    \put(150,170){\circle*{6}}
    \put(180,170){\circle*{6}}
    \put(60,195){\circle*{6}}
    \put(0,170){\line(1,0){180}}
    \put(60,170){\line(0,1){25}}
    \put(0,165){\makebox(0,0)[t]{\scriptsize$a_1$}}
    \put(30,165){\makebox(0,0)[t]{\scriptsize$a_2$}}
    \put(60,165){\makebox(0,0)[t]{\scriptsize$a_3$}}
    \put(90,165){\makebox(0,0)[t]{\scriptsize$a_4$}}
    \put(120,165){\makebox(0,0)[t]{\scriptsize$a_5$}}
    \put(150,165){\makebox(0,0)[t]{\scriptsize$a_6$}}
    \put(180,165){\makebox(0,0)[t]{\scriptsize$a_7$}}
    \put(70,195){\makebox(0,0)[c]{\scriptsize$a_8$}}

    \put(-50,130){\makebox(0,0)[c]{$(F_4)_{odd}$}}
    \put(0,130){\circle*{6}}
    \put(30,130){\circle*{6}}
    \put(60,130){\circle*{6}}
    \put(90,130){\circle*{6}}
    \put(0,130){\line(1,0){30}}
    \put(60,130){\line(1,0){30}}
    \put(0,125){\makebox(0,0)[t]{\scriptsize$a_1$}}
    \put(30,125){\makebox(0,0)[t]{\scriptsize$a_2$}}
    \put(60,125){\makebox(0,0)[t]{\scriptsize$a_3$}}
    \put(90,125){\makebox(0,0)[t]{\scriptsize$a_4$}}

    \put(-50,90){\makebox(0,0)[c]{$(H_3)_{odd}$}}
    \put(0,90){\circle*{6}}
    \put(30,90){\circle*{6}}
    \put(60,90){\circle*{6}}
    \put(0,90){\line(1,0){60}}
    \put(15,95){\makebox(0,0)[b]{$5$}}
       \put(0,85){\makebox(0,0)[t]{\scriptsize$a_1$}}
       \put(30,85){\makebox(0,0)[t]{\scriptsize$a_2$}}
       \put(60,85){\makebox(0,0)[t]{\scriptsize$a_3$}}

    \put(-50,50){\makebox(0,0)[c]{$(H_4)_{odd}$}}
    \put(0,50){\circle*{6}}
    \put(30,50){\circle*{6}}
    \put(60,50){\circle*{6}}
    \put(90,50){\circle*{6}}
    \put(0,50){\line(1,0){90}}
    \put(15,55){\makebox(0,0)[b]{$5$}}
       \put(0,45){\makebox(0,0)[t]{\scriptsize$a_1$}}
       \put(30,45){\makebox(0,0)[t]{\scriptsize$a_2$}}
       \put(60,45){\makebox(0,0)[t]{\scriptsize$a_3$}}
       \put(90,45){\makebox(0,0)[t]{\scriptsize$a_4$}}

    \put(-50,15){\makebox(0,0)[c]{$(I_2(m))_{odd}~~(m \geq 5)$}}
    \put(20,15){\circle*{6}}
    \put(50,15){\circle*{6}}
    \put(20,15){\line(1,0){30}}
    \put(35,20){\makebox(0,0)[b]{$m$}}
       \put(20,10){\makebox(0,0)[t]{\scriptsize$a_1$}}
       \put(50,10){\makebox(0,0)[t]{\scriptsize$a_2$}}
    \put(100,15){\makebox(0,0)[c]{$m$ odd}}


    \put(20,-9){\circle*{6}}
    \put(50,-9){\circle*{6}}
       \put(20,-14){\makebox(0,0)[t]{\scriptsize$a_1$}}
       \put(50,-14){\makebox(0,0)[t]{\scriptsize$a_2$}}
    \put(100,-9){\makebox(0,0)[c]{$m$ even}}

\end{picture}
  \caption{$\Gamma_{odd}$ for the irreducible spherical-type Coxeter graphs $\Gamma$}
  \label{fig:odd-coxeter-graphs}
 \end{figure}

Figure \ref{fig:odd-coxeter-graphs} shows that each irreducible
spherical-type Artin group falls into one of two classes; (i) those
in which $\Gam_{odd}$ is connected and (ii) those in which
$\Gam_{odd}$ has two components.  Within a given class the
arguments are quite similar. Thus, we will only show the complete
details of the computations for types $A_n$ and $B_n$.  The rest
of the types have similar computations.

\subsubsection{Lemmas for simplifying presentations}
\label{subsec:two-lemmas}

 We will encounter two sets of relations quite often in our
 computations and it will be necessary to replace them with sets
 of simpler but equivalent relations.  In this section we give two
 lemmas which allow us to make these replacements.

 Let $\{ p_k \}_{k\in \Z}$ , $a$, $b$, and $q$ be letters.  In the
 following keep in mind that the relators
 $p_{k+1}p_{k+2}^{-1}p_{k}^{-1}$ split up into the two types of
 relations $p_{k+2}=p_{k}^{-1} p_{k+1}$ (for $k\geq0$), and
 $p_{k}=p_{k+1}p_{k+2}^{-1}$ (for $k<0$).  The two lemmas are:

 \begin{lemma}
 \label{lem:relation-replacement-1}
   The set of relations
    \begin{eqnarray}
      p_{k+1}p_{k+2}^{-1}p_{k}^{-1}=1, \quad
      p_k a p_{k+2} a^{-1} p_{k+1}^{-1} a^{-1}=1, \quad
      b=p_0ap_{0}^{-1},
      \label{relations-set-1}
    \end{eqnarray}
   is equivalent to the set
    \begin{eqnarray}
      p_{k+1}p_{k+2}^{-1}p_{k}^{-1}&=&1, \label{relations-set-2-1}\\
      p_0ap_{0}^{-1}&=& b, \label{relations-set-2-2}\\
      p_0bp_{0}^{-1}&=& b^2a^{-1}b \label{relations-set-2-3} \\
      p_1ap_{1}^{-1}&=& a^{-1}b,  \label{relations-set-2-4} \\
      p_1bp_{1}^{-1}&=& (a^{-1}b)^{3}a^{-2}b. \label{relations-set-2-5}
    \end{eqnarray}
\end{lemma}

\begin{lemma}
\label{lem:relation-replacement-2}
  The set of relations:
    \begin{center}
       $p_{k+1}p_{k+2}^{-1}p_{k}^{-1}=1$, \quad $p_kq=qp_{k+1}$,
    \end{center}
    is equivalent to the set
    \begin{center}
       $p_{k+1}p_{k+2}^{-1}p_{k}^{-1}=1$,\quad $p_0q=qp_1$,\quad
       $p_1q = qp_{0}^{-1}p_1$.
    \end{center}
\end{lemma}

The proof of lemma~\ref{lem:relation-replacement-2} is
straightforward.  On the other hand, the proof of the
lemma~\ref{lem:relation-replacement-1} is somewhat long and
tedious.

\begin{proof}[Lemma~\ref{lem:relation-replacement-2}]
  Clearly the second set of relations follows from the first set of
  relations since $p_2=p_{0}^{-1}p_{1}$.  To prove
  the converse we first prove that $p_kq=qp_{k+1}$ ($k\geq 0$)
  follows from the second set of relations by induction on $k$.
  It is easy to see then that the same is true for $k<0$.  For
  $k=0,1$ the result clearly holds.  Now, for $k=m+2$;
  \begin{eqnarray*}
   p_{m+2}qp_{m+3}^{-1}q^{-1}&=&p_{m+2}qp_{m+2}^{-1}p_{m+1}q^{-1},\\
                                   &=&p_{m+2}(p_{m+1}^{-1}q)p_{m+1}q^{-1} \text{~~~by IH $(k=m+1)$},\\
                                   &=&p_{m+2}p_{m+1}^{-1}(qp_{m+1})q^{-1},\\
                                   &=&p_{m+2}p_{m+1}^{-1}(p_mq)q^{-1} \text{~~~by IH $(k=m)$},\\
                                   &=&p_{m+2}p_{m+1}^{-1}p_{m},\\
                                   &=&1.
  \end{eqnarray*}
\end{proof}

\begin{proof}[Lemma~\ref{lem:relation-replacement-1}]
  First we show the second set of relations follows from the
  first set.
  Taking $k=0$ in the second relation in (\ref{relations-set-1})
  we get the relation
  \begin{eqnarray*}
   p_0 a p_2 a^{-1} p_{1}^{-1}a^{-1}=1,
  \end{eqnarray*}
  and, using the relations $p_2=p_{0}^{-1}p_1$ and
  $b=p_0ap_{0}^{-1}$, (\ref{relations-set-2-4}) easily follows.
  Taking $k=1$ in the second relation in (\ref{relations-set-1}) we
  get the relation
  \begin{eqnarray*}
   p_1 a p_3 a^{-1} p_{2}^{-1}a^{-1}=1.
  \end{eqnarray*}
  Using the relations $p_3=p_{1}^{-1}p_2$ and $p_2=p_{0}^{-1}p_1$ this becomes
  \begin{eqnarray*}
   p_1 a p_{1}^{-1}p_{0}^{-1}p_1 a^{-1} p_{1}^{-1}p_0a^{-1}=1.
  \end{eqnarray*}
  But $p_1ap_{1}^{-1}=a^{-1}b$ (by (\ref{relations-set-2-4})) so this reduces to
  \begin{eqnarray*}
   a^{-1}b p_{0}^{-1}b^{-1}a p_0a^{-1}=1.
  \end{eqnarray*}
  Isolating $bp_{0}^{-1}$ on one side of the equation gives
  \begin{eqnarray*}
    bp_{0}^{-1}=a^{2}p_{0}^{-1}a^{-1}b.
  \end{eqnarray*}
  Multiplying both sides on the left by $p_0$ and using the relation
  $p_0ap_{0}^{-1}=b$ it easily follows $p_0bp_{0}^{-1}=b^{2}a^{-1}b$, which is
  (\ref{relations-set-2-3}).
  Finally, taking $k=2$ in the second relation in (\ref{relations-set-1}) we get the relation
  \begin{eqnarray*}
   p_2 a p_4 a^{-1} p_{3}^{-1}a^{-1}=1.
  \end{eqnarray*}
  Using the relation $p_4=p_{2}^{-1}p_{3}$ this becomes
  \begin{eqnarray}
  \label{transitional-relation}
   p_2 a p_{2}^{-1}p_3 a^{-1} p_{3}^{-1}a^{-1}=1.
  \end{eqnarray}
  Note that
  \begin{eqnarray*}
   p_2ap_{2}^{-1}&=&p_{0}^{-1}p_1ap_{1}^{-1}p_0 \text{\quad by $p_2=p_{0}^{-1}p_{1}$}\\
                 &=&p_{0}^{-1}a^{-1}bp_{0} \text{\quad by (\ref{relations-set-2-4})}\\
                 &=&a^{-2}ba^{-1}a \text{\quad using (\ref{relations-set-1}) and (\ref{relations-set-2-3})}\\
                 &=&a^{-2}b
  \end{eqnarray*}
  and
  \begin{eqnarray*}
   p_3ap_{3}^{-1}&=&p_{1}^{-1}p_{2}ap_{2}^{-1}p_1 {\rm~~~~~by~p_3=p_{1}^{-1}p_{2}}\\
                 &=&p_{1}^{-1}a^{-2}bp_1,
  \end{eqnarray*}
  where the second equality follows from the previous statement.
  Thus, (\ref{transitional-relation}) becomes
  \begin{eqnarray*}
   a^{-2}bp_{1}^{-1}b^{-1}a^{2}p_1a^{-1}=1
  \end{eqnarray*}
  Isolating the factor $bp_{1}^{-1}$ on one side of the equation,
  multiplying both sides by $p_1$, and using the
  relation (\ref{relations-set-2-4}) we easily get the relation
  (\ref{relations-set-2-5}).  Therefore we have
  that the second set of relations
  (\ref{relations-set-2-1})-(\ref{relations-set-2-5})
  follows from the first set of relations (\ref{relations-set-1}).

  In order to show the relations in (\ref{relations-set-1}) follow
  from the relations in
  (\ref{relations-set-2-1})-(\ref{relations-set-2-5}) it suffices
  to just show that the second relation in (\ref{relations-set-1})
  follows from the relations in
  (\ref{relations-set-2-1})-(\ref{relations-set-2-5}).  To do this
  we need the following fact: The relations
  \begin{eqnarray}
     p_k a p_{k}^{-1}&=& a^{k}b, \label{relations-set-3-1} \\
     p_k b p_{k}^{-1}&=& (a^{-k}b)^{k+2}a^{-(k+1)}b, \label{relations-set-3-2} \\
     p_{k}^{-1} a p_{k}&=& ab^{-1}a^{k+2},  \label{relations-set-3-3} \\
     p_{k}^{-1} b p_{k}&=& (ab^{-1}a^{k+2})^k a,
     \label{relations-set-3-4}
    \end{eqnarray}
  follow from the relations in
  (\ref{relations-set-2-1})-(\ref{relations-set-2-5}).
  The proof of this fact is left to
  lemma ~\ref{lem:relation-replacement-3} below.
  From the relations (\ref{relations-set-3-1})-(\ref{relations-set-3-4}) we obtain
  \begin{eqnarray}
  \label{rel-1}
    p_{k+1}ap_{k+1}^{-1}=a^{-(k+1)}b=a^{-1}\cdot a^{-k}b =
    a^{-1}p_{k}ap_{k}^{-1},
  \end{eqnarray}
  and
  \begin{eqnarray}
  \label{rel-2}
    p_{k+1}^{-1}ap_{k+1}=ab^{-1}a^{k+3}=ab^{-1}a^{k+2}a
    =p_{k}^{-1}ap_ka.
  \end{eqnarray}
  Now we are in a position to show that that the second relation in
  (\ref{relations-set-1}) follows from the relations in
  (\ref{relations-set-2-1})-(\ref{relations-set-2-5}).  For $k \geq 0$
  \begin{eqnarray*}
  p_kap_{k+2}a^{-1}p_{k+1}^{-1}a^{-1}&=&p_kap_{k}^{-1} \underbrace{p_{k+1}a^{-1}p_{k+1}^{-1}}a^{-1}
                                          \text{~~~~~by (\ref{relations-set-2-1})}\\
                                     &=&p_kap_{k}^{-1}(a^{-1}p_kap_{k}^{-1})^{-1}a^{-1}
                                          \text{~~~~~by (\ref{rel-1})}\\
                                     &=&1.
  \end{eqnarray*}
  and for $k<0$
  \begin{eqnarray*}
  p_kap_{k+2}a^{-1}p_{k+1}^{-1}a^{-1}&=&p_{k+1}\underbrace{p_{k+2}^{-1}ap_{k+2}}a^{-1}p_{k+1}^{-1}a^{-1}
                                         \text{~~~~~by (\ref{relations-set-2-1})}  \\
                                     &=&p_{k+1}(p_{k+1}^{-1}ap_{k+1}a)a^{-1}p_{k+1}^{-1}a^{-1}
                                         \text{~~~~~by (\ref{rel-2})} \\
                                     &=&1.
  \end{eqnarray*}
  Therefore, the relations
  \begin{eqnarray*}
    p_kap_{k+2}a^{-1}p_{k+1}^{-1}a^{-1}=1 ,~~~~ k \in \Z
  \end{eqnarray*}
  follow from the relations in
  (\ref{relations-set-2-1})-(\ref{relations-set-2-5}).
\end{proof}

To complete the proof of lemma~\ref{lem:relation-replacement-1} we
need to prove the following.

\begin{lemma}
\label{lem:relation-replacement-3}
  The relations
    \begin{eqnarray*}
     p_k a p_{k}^{-1}&=& a^{k}b \\
     p_k b p_{k}^{-1}&=& (a^{-k}b)^{k+2}a^{-(k+1)}b \\
     p_{k}^{-1} a p_{k}&=& ab^{-1}a^{k+2} \\
     p_{k}^{-1} b p_{k}&=& (ab^{-1}a^{k+2})^k a
    \end{eqnarray*}
   follow from the relations in
   (\ref{relations-set-2-1})-(\ref{relations-set-2-5}).
\end{lemma}
  \begin{proof}
    We will use induction to prove the result for nonnegative
    indices k, the result for negative indices k is similar.
    Clearly this holds for $k=0,1$. For $k=m+2$ we have
    \begin{eqnarray*}
      p_{m+2}ap_{m+2}^{-1}&=&p_{m}^{-1}p_{m+1}ap_{m+1}^{-1}p_{m}
      \text{~~~~~by (\ref{relations-set-2-1}),}\\
                          &=&p_{m}^{-1}a^{-(m+1)}bp_{m} \text{~~~~~by induction hypothesis (IH),}\\
                          &=&(p_{m}^{-1}a^{-(m+1)}p_m)(p_{m}^{-1}bp_{m}),\\
                          &=&(p_{m}^{-1}ap_m)^{-(m+1)}(p_{m}^{-1}bp_{m}),\\
                          &=&(ab^{-1}a^{m+2})^{-(m+1)}(ab^{-1}a^{m+2})^{m}a \text{~~~~~by IH,}\\
                          &=&(ab^{-1}a^{m+2})^{-1}a,\\
                          &=&a^{-(m+2)}b,
    \end{eqnarray*}
    \begin{eqnarray*}
      p_{m+2}bp_{m+2}^{-1}&=&p_{m}^{-1}p_{m+1}bp_{m+1}^{-1}p_{m}
      \text{~~~~~by (\ref{relations-set-2-1}),}\\
                          &=&p_{m}^{-1}(a^{-(m+1)}b)^{m+3}a^{-(m+2)}bp_{m} \text{~~~~~by IH,}\\
                          &=&((p_{m}^{-1}ap_m)^{-(m+1)}(p_{m}^{-1}bp_m))^{m+3}(p_{m}^{-1}ap_{m})^{-(m+2)}p_{m}^{-1}bp_{m},\\
                          &=&((ab^{-1}a^{m+2})^{-(m+1)}(ab^{-1}a^{m+2})^{m}a)^{(m+3)} \\
                          &~& ~~~~ \cdot (ab^{-1}a^{m+2})^{-(m+2)}(ab^{-1}a^{m+2})^m a \text{~~~~~by IH,}\\
                          &=&(a^{-(m+2)}b)^{m+3}(ab^{-1}a^{m+2})^{-2}a,\\
                          &=&(a^{-(m+2)}b)^{m+4}a^{-(m+3)}b,
    \end{eqnarray*}
  Similarly for the other two equations.
    Thus, the result follows by induction.
  \end{proof}

\subsubsection{Type $A$}
\label{subsec:type-A}

The first presentation for the commutator subgroup $\Br_{n+1}'=
\A_{A_n}'$ of the braid group $\Br_{n+1} =\A_{A_n}$ appeared in
\cite{GL69} but the details of the computation were minimal. Here
we fill in the details of Gorin and Lin's computation.

The presentation of $\A_{A_n}$ is
 \begin{eqnarray*}
    \A_{A_n}=\la a_1,...,a_{n}: &~&a_{i}a_{j}=a_{j}a_{i}
    \text{~~~for $|i-j|\geq 2$},\\
    &~&a_{i}a_{i+1}a_{i}=a_{i+1}a_{i}a_{i+1}
    \text{~~~for $1 \leq i \leq n-1$ } \ra.
 \end{eqnarray*}
Since $(A_n)_{odd}$ is connected then by proposition
\ref{thm:characterization-commuator-subgroup} $\A_{A_n}' = \ker
(\d$). Elements $U,V \in \A_{A_n}$ lie in the same right coset of
$\A_{A_n}'$ precisely when they have the same degree:
  \begin{eqnarray*}
      \A_{A_n}' U = \A_{A_n}' V &\Longleftrightarrow& U
      V^{-1} \in \A_{A_n}' \\ &\Longleftrightarrow&
      \d(U)=\d(V),
  \end{eqnarray*}
thus a Schreier system of right coset representatives for
$\A_{A_n}$ modulo $\A_{A_n}'$ is
  \begin{eqnarray*}
    \RR = \{a_{1}^{k} : k \in \mathbb Z\}
  \end{eqnarray*}
By the Reidemeister-Schreier method, in particular equation
(\ref{general-presentation}), $\A_{A_n}'$ has generators
$s_{a_{1}^{k},a_{j}}:=
a_{1}^{k}a_{j}(\overline{a_{1}^{k}a_{j}})^{-1}$ with presentation
\begin{eqnarray}
\label{A-prelim-presentation}
   \la s_{a_{1}^{k},a_{j}}, \ldots : s_{a_{1}^{m},a_{\lambda}} , \ldots, \tau(a_{1}^{\ell}R_{i}a_{1}^{-\ell}),\ldots,
   \tau(a_{1}^{\ell}T_{i,j}a_{1}^{-\ell}), \ldots \ra ,
\end{eqnarray}
 where $j\in\{1, \ldots,n\},$ $k,\ell \in \Z$, and
 $m \in \Z$, $\lambda \in \{1, \ldots,n\}$ such that
 $a_{1}^{m} a_{\lambda} \thickapprox \overline{a_{1}^{m}
 a_{\lambda}}$ (``freely equal"), and $T_{i,j}$, $R_{i}$ represent the
 relators
 $a_{i}a_{j}a_{i}^{-1}a_{j}^{-1}$, $|i-j| \geq 2$, and
 $a_{i} a_{i+1} a_{i} a_{i+1}^{-1} a_{i}^{-1} a_{i+1}^{-1}$,
 respectively.  Our goal is to clean up this presentation.

 The first thing to notice is that
  \begin{eqnarray*}
   a_{1}^{m} a_{\lambda} \thickapprox \overline{a_{1}^{m}
   a_{\lambda}} = a_{1}^{m+1} \Longleftrightarrow \lambda =1
  \end{eqnarray*}
 Thus, the first type of relation
 in (\ref{A-prelim-presentation}) is precisly $s_{a_{1}^{m},a_{1}}=1$, for all $m \in \Z$.

 Next, we use the definition of the Reidemeister rewriting
 function (\ref{reidemeister-rewriting-function}) to express the second and third types of
 relations in (\ref{A-prelim-presentation}) in terms of the generators
 $s_{a_{1}^{k},a_{j}}$:
  \begin{eqnarray}
    \tau(a_{1}^{k}T_{i,j}a_{1}^{-k})&=&s_{a_{1}^{k},a_{i}}s_{a_{1}^{k+1},a_{j}}s_{a_{1}^{k+1},a_{i}}^{-1}s_{a_{1}^{k},a_{j}}^{-1}  \label{A-prelim-relations-1} \\
    \tau(a_{1}^{k}R_{i}a_{1}^{-k})&=&s_{a_{1}^{k},a_{i}}s_{a_{1}^{k+1},a_{i+1}}s_{a_{1}^{k+2},a_{i}}s_{a_{1}^{k+2},a_{i+1}}^{-1}s_{a_{1}^{k+1},a_{i}}^{-1}s_{a_{1}^{k},a_{i+1}}^{-1}
    \label{A-prelim-relations-2}
  \end{eqnarray}
 Taking $i=1$, $j \geq 3$ in (\ref{A-prelim-relations-1}) we get
  \begin{eqnarray*}
    s_{a_{1}^{k+1},a_{j}}=s_{a_{1}^{k},a_{j}}
  \end{eqnarray*}
 Thus, by induction on $k$,
  \begin{eqnarray}
  \label{A-q-relations}
    s_{a_{1}^{k},a_{j}}=s_{1,a_{j}}
  \end{eqnarray}
 for $j \geq 3$ and for all $k \in \Z$.

 Therefore, $\A_{A_n}'$ is generated by
 $s_{a_{1}^{k},a_{2}}=a_{1}^{k}a_{2}a_{1}^{-(k+1)}$
 and $s_{1,a_{\ell}}= a_{\ell}a_{1}^{-1}$, where $k \in \Z$ , $3
 \leq \ell \leq n$.  To simplify notation let us rename the
 generators; let $p_k:=a_{1}^{k}a_{2}a_{1}^{-(k+1)}$ and
 $q_{\ell}:=a_{\ell}a_{1}^{-1}$, for $k \in \Z$ , $3
 \leq \ell \leq n$.  We now investigate the relations in
 (\ref{A-prelim-relations-1}) and (\ref{A-prelim-relations-2}).

 The relations in (\ref{A-prelim-relations-2}) break up into the following
 three types (using ~\ref{A-q-relations}):
  \begin{eqnarray}
    p_{k+1} p_{k+2}^{-1} p_{k}^{-1} &~& \text{(taking $i=1$)}
    \label{A-relation-type-1-1}\\
    p_k q_3 p_{k+2} q_{3}^{-1} p_{k+1}^{-1} q_{3}^{-1} &~& \text{(taking $i=2$)}
    \label{A-relation-type-1-2}\\
    q_i q_{i+1} q_{i} q_{i+1}^{-1} q_{i}^{-1} q_{i+1}^{-1} &~& \text{for $3\leq i \leq
    n-1$.}
    \label{A-relation-type-1-3}
  \end{eqnarray}
 The relations in (~\ref{A-prelim-relations-1}) break up into the following two
 types:
  \begin{eqnarray}
     p_k q_j p_{k+1}^{-1} q_{j}^{-1}&~& \text{for $4 \leq j \leq n$, (\text{taking }$i=2$)}
     \label{A-relation-type-2-1}\\
     q_i q_j q_{i}^{-1} q_{j}^{-1} &~& \text{for $3 \leq i<j \leq n$, $|i-j| \geq 2$.}
     \label{A-relation-type-2-2}
  \end{eqnarray}

 We now have a presentation for $\A_{A_n}'$ consisting of the
 generators $p_k,q_{\ell}$, where $k \in \Z$, $3
 \leq \ell \leq n-1$, and defining relations (\ref{A-relation-type-1-1})
 -(\ref{A-relation-type-2-2}).  However, notice
 that relation (\ref{A-relation-type-1-1}) splits up into the two relations
  \begin{eqnarray}
    p_{k+2}=p_{k}^{-1} p_{k+1} &~& \text{for $k\geq0$,}\\
    p_{k}=p_{k+1}p_{k+2}^{-1}  &~& \text{for $k<0$.}
  \end{eqnarray}
  Thus, for $k \neq 0,1,$ $p_k$ can be expressed in terms of $p_0$
  and $p_1$.  It follows that $\A_{A_n}'$ is finitely
  generated.  In order to show $\A_{A_n}'$ is finitely presented we
  need to be able to replace the infinitly many relations in (\ref{A-relation-type-1-2})
  and (\ref{A-relation-type-2-1}) with finitely many relations.
  This can be done using lemmas \ref{lem:relation-replacement-1}
  and \ref{lem:relation-replacement-2}, but this requires us to add a new letter $b$
  to the generating set with a new relation $b=p_0q_3p_{0}^{-1}$.
  Thus $\A_{A_n}'$ is generated by $p_0,p_1,q_{\ell},b$, where $3\leq \ell \leq n-1$,
  with defining relations:
  \begin{gather*}
    p_0q_3p_{0}^{-1}= b, \quad p_0bp_{0}^{-1}= b^2q_3^{-1}b,
    \quad
    p_1q_3p_{1}^{-1}= q_3^{-1}b, \quad p_1bp_{1}^{-1}= (q_3^{-1}b)^{3}q_3^{-2}b, \\
    q_i q_{i+1} q_{i} q_{i+1}^{-1} q_{i}^{-1} q_{i+1}^{-1} \text{~~~~($3\le i \le n-1$)},\\
    p_0q_j=q_jp_1 \text{~~~~($4\le j \le n$)}, \quad p_1q_j= q_jp_{0}^{-1}p_1 \text{~~~~($4\le j \le n$)}. \\
    q_i q_j q_{i}^{-1} q_{j}^{-1} \text{~~~~~~~~($3 \leq i<j \leq n$, $|i-j| \geq 2$).}
  \end{gather*}

 Noticing that for $n=2$ the generators $q_k$ ($3\le k \le n$),
 and $b$ do not exist, and for $n=3$ the generators $q_k$ ($4 \le
 k \le n$) do not exist, we have proved the following theorem.

 \begin{theorem}
 \label{thm:commutator-type-A-presentation}
   For every $n\geq2$ the commutator subgroup $\A_{A_n}'$ of the
   Artin group $\A_{A_n}$ is a finitely  presented group.  $\A_{A_2}'$ is a free
   group with two free generators
   \begin{eqnarray*}
   p_0=a_{2}a_{1}^{-1},~~
   p_1=a_{1}a_{2}a_{1}^{-2}.
   \end{eqnarray*}
   $\A_{A_3}'$ is the group generated by
   \begin{eqnarray*}
    p_0=a_{2}a_{1}^{-1},~~
    p_1=a_{1}a_{2}a_{1}^{-2},~~
    q=a_{3}a_{1}^{-1},~~
    b=a_{2}a_{1}^{-1}a_{3}a_{2}^{-1},
   \end{eqnarray*}
   with defining relations
   \begin{gather*}
    b=p_0qp_{0}^{-1}, \quad
   p_0bp_{0}^{-1}=b^{2}q^{-1}b, \\
   p_1qp_{1}^{-1}=q^{-1}b, \quad
   p_1bp_{1}^{-1}=(q^{-1}b)^{3}q^{-2}b.
   \end{gather*}
   For $n\geq4$ the group $\A_{A_n}'$ is generated by
    \begin{eqnarray*}
    p_0=a_{2}a_{1}^{-1},~~
    p_1=a_{1}a_{2}a_{1}^{-2},~~
    q_3=a_{3}a_{1}^{-1},~~ \\
    b=a_{2}a_{1}^{-1}a_{3}a_{2}^{-1},~~
    q_{\ell}=a_{\ell}a_{1}^{-1} ~~(4\leq \ell \leq n-1),
   \end{eqnarray*}
   with defining relations
   \begin{gather*}
    b=p_0q_3p_{0}^{-1}, \quad
    p_0bp_{0}^{-1}=b^{2}q_3^{-1}b,\\
    p_1q_3p_{1}^{-1}=q_3^{-1}b, \quad
    p_1bp_{1}^{-1}=(q_3^{-1}b)^{3}q_3^{-2}b,\\
    p_0q_i=q_ip_1 ~~(4\leq i \leq n), \quad p_1q_i=q_ip_{0}^{-1}p_{1} ~~(4\leq i \leq n)\\
    q_3q_i=q_iq_3 ~~(5\leq i \leq n), \quad q_3q_4q_3=q_4q_3q_4,\\
    q_iq_j=q_jq_i ~~(4 \leq i< j-1\leq n-1), \quad q_iq_{i+1}q_i=q_{i+1}q_{i}q_{i+1}~~(4 \leq i\leq n-1).
   \end{gather*}
   \hfill $\Box$
  \end{theorem}

  \begin{corollary}
  \label{cor:commutator-type-A-perfect}
    For $n\ge4$ the commutator subgroup $\A_{A_n}'$ of the Artin
    group of type $A_n$ is finitely generated and perfect (i.e. $\A_{A_n}''=\A_{A_n}'$).
  \end{corollary}

  \begin{proof}
    Abelianizing the presentation of $\A_{A_n}'$ in the theorem
    results in a presentation of the trivial group.  Hence
    $\A_{A_n}''=\A_{A_n}'$.
  \end{proof}

 Now we study in greater detail the group $\A_{A_3}'$, the results
 of which will be used in section \ref{subsec:type-A-indicability}.
 From the presentation of $\A_{A_3}'$ given in theorem
 \ref{thm:commutator-type-A-presentation} one can easily deduce
 the relations:
 \begin{alignat*}{2}
   p_{0}^{-1}qp_{0} & =qb^{-1}q^{2}, & \quad
   p_{0}^{-1}bp_{0} & =q, \\
   p_{1}^{-1}qp_{1} & =qb^{-1}q^{3}, & \quad
   p_{1}^{-1}bp_{1} & =qb^{-1}q^{4}.
 \end{alignat*}
 Let $T$ be the subgroup of $\A_{A_3}'$ generated by $q$ and $b$.
 The above relations and the defining relations in the
 presentation for $\A_{A_3}'$ tell us that $T$ is a normal
 subgroup of $\A_{A_3}'$.  To obtain a representation of the
 factor group $\A_{A_3}'/T$ it is sufficient to add to the
 defining relations in the presentation for $\A_{A_3}'$ the
 relations $q=1$ and $b=1$.  It is easy to see this results in the
 presentation of the free group generated by $p_0$ and $p_1$.
 Thus, $\A_{A_3}'/T$ is a free group of rank $2$, $F_2$.  We
 have the exact sequence
 $$
 \begin{CD}
  1 \longrightarrow T \longrightarrow \A_{A_3}' \longrightarrow
  \A_{A_3}'/T \longrightarrow 1.
 \end{CD}
 $$
 Since $\A_{A_3}'/T$ is free then the exact sequence is actually
 split so
 \begin{eqnarray*}
  \A_{A_3}' \simeq T \rtimes \A_{A_3}'/T  \simeq T \rtimes F_2,
 \end{eqnarray*}
 where the action of $F_2$ on $T$ is given by the defining
 relations in the presentation of $\A_{A_3}'$ and the relations
 above.  In \cite{GL69} it is shown (theorem 2.6) the group $T$ is also free of rank 2,
 so we have the following theorem.

 \begin{theorem}
 \label{thm:commutator-subgroup-A3}
  The commutator subgroup $\A_{A_3}'$ of the Artin group of type
  $A_3$ is the semidirect product of two free groups each of rank 2;
  \begin{eqnarray*}
    \A_{A_3}' \simeq F_2 \rtimes F_2.
  \end{eqnarray*}
  \hfill $\Box$
 \end{theorem}

\subsubsection{Type $B$}
\label{subsec:type-B}

The presentation of $\A_{B_n}$ is
 \begin{eqnarray*}
    \A_{B_n}=\la a_1,...,a_{n}: &~&a_{i}a_{j}=a_{j}a_{i}
    \text{~~~for $|i-j|\geq 2$},\\
    &~&a_{i}a_{i+1}a_{i}=a_{i+1}a_{i}a_{i+1}
    \text{~~~for $1 \leq i \leq n-2$ } \\
    &~&a_{n-1}a_{n}a_{n-1}a_{n}=a_{n}a_{n-1}a_{n}a_{n-1}
    \ra.
 \end{eqnarray*}
 Let $T_{i,j},R_i$ ($1\le i \le n-2$), and $R_{n-1}$ denote the
 associated relators $a_{i} a_{j} a_{i}^{-1} a_{j}^{-1}$,
 $a_{i} a_{i+1} a_{i} a_{i+1}^{-1} a_{i}^{-1} a_{i+1}^{-1}$, and
 $a_{n-1} a_{n} a_{n-1} a_{n} a_{n-1}^{-1} a_{n}^{-1} a_{n-1}^{-1} a_{n}^{-1}$,
 respectively.

 As seen in figure \ref{fig:odd-coxeter-graphs} the graph
 $(B_n)_{odd}$ has two components: $\Gam_1$ and
 $\Gam_2$, where $\Gam_2$ denotes the component containing the
 single vertiex $a_n$.  Let $\d_1$ and $\d_2$ denote the
 associated degree maps, respectively, so from proposition ~\ref{thm:characterization-commuator-subgroup}
 \begin{eqnarray*}
 \A_{B_n}'=\{U\in \A_{B_n} : \d_1(U)=0 \text{~and~} \d_2(U)=0
 \}.
 \end{eqnarray*}
 For elements $U,V \in \A_{A_n}$,
  \begin{eqnarray*}
      \A_{B_n}' U = \A_{B_n}' V &\Leftrightarrow& U
      V^{-1} \in \A_{B_n}' \\ &\Leftrightarrow&
      \d_1(U)=\d_1(V), \text{~and}\\
      &~&\d_2(U)=\d_2(V),
  \end{eqnarray*}
 thus a Schreier system of right coset representatives for
 $\A_{B_n}$ modulo $\A_{B_n}'$ is
  \begin{eqnarray*}
    \RR = \{a_{1}^{k}a_{n}^{\ell} : k,\ell \in \mathbb Z \}
  \end{eqnarray*}
 By the Reidemeister-Schreier method, in particular equation
 (\ref{general-presentation}), $\A_{B_n}'$ is generated by
 \begin{eqnarray*}
  s_{a_{1}^{k}a_{n}^{k},a_{j}}
      &:=& a_{1}^{k}a_{n}^{\ell}a_{j}(\overline{a_{1}^{k}a_{n}^{\ell}a_{j}})^{-1}\\
      &=&
      \begin{cases}
       a_{1}^{k}a_{n}^{\ell}a_{j}a_{n}^{-\ell}a_{1}^{-(k+1)}  \text{~~~if $j\ne n$}\\
       1  \text{~~~~~~~~~~~~~~~~~~~~~~~~~if $j=n$.}
      \end{cases}
 \end{eqnarray*}
 with presentation
 \begin{equation}
 \label{B-prelim-presentation}
   \begin{split}
     \A_{B_n}' = \la s_{a_{1}^{k}a_{n}^{\ell},a_{j}}, \ldots : &~ s_{a_{1}^{p}a_{n}^{q},a_{\lambda}} ,
     \ldots,\\
     &~ \tau(a_{1}^{k}a_{n}^{\ell}T_{i,j}(a_{1}^{k}a_{n}^{\ell})^{-1}),\ldots, \text{~~~($1\le i<j\le n,|i-j|\ge2$)},\\
     &~ \tau(a_{1}^{k}a_{n}^{\ell}R_{i}(a_{1}^{k}a_{n}^{\ell})^{-1}),\ldots, \text{~~~($1\le i \le n-2$)},\\
     &~ \tau(a_{1}^{k}a_{n}^{\ell}R_{n-1}(a_{1}^{k}a_{n}^{\ell})^{-1}),\ldots \ra ,
   \end{split}
 \end{equation}
 where $p,q \in \Z$, $\lambda \in \{1, \ldots,n-1\}$ such that
 $a_{1}^{p}a_{n}^{q} a_{\lambda} \thickapprox
 \overline{a_{1}^{p}a_{n}^{q} a_{\lambda}}$ ("freely equal").
 Again, our goal is to clean up this presentation.

 The cases $n=2,3$, and $4$ are straightforward
 after one sees the computation for the general case $n\ge 5$, so
 we will not include the computations for these cases.  The results are included in
 theorem \ref{thm:commutator-type-B-presentation}.  From now on it
 will be assumed that $n\ge5$.

 Since
 \begin{equation*}
   a_{1}^{p}a_{n}^{q} a_{\lambda} \thickapprox
   \overline{a_{1}^{p}a_{n}^{q} a_{\lambda}} =\\
   \begin{cases}
    a_{1}^{p+1}a_{n}^{q}   \text{~~~$\lambda\ne n$} \\
    a_{1}^{p}a_{n}^{q+1}   \text{~~~$\lambda=n$}
   \end{cases}
   \Longleftrightarrow
   \lambda=n \text{~or;~} \\
   \lambda=1 \text{~and~} q=0,
 \end{equation*}
 the first type of
 relations in (\ref{B-prelim-presentation}) are precisely
 \begin{eqnarray}
 \label{B-relation-type-1}
    s_{a_{1}^{k}a_{n}^{\ell},a_{n}}=1, \text{~~and~~}
    s_{a_{1}^{k},a_{1}}=1.
 \end{eqnarray}

 The second type of relations in (\ref{B-prelim-presentation}), after rewriting using
 equation (\ref{reidemeister-rewriting-function}), are
 \begin{equation}
 \label{B-relation-type-2}
   s_{a_{1}^{k}a_{n}^{\ell},a_{i}}s_{\overline{a_{1}^{k}a_{n}^{\ell}a_{i}},a_{j}}
   s_{\overline{a_{1}^{k}a_{n}^{\ell}a_{i}a_{j}a_{i}^{-1}},a_{i}}^{-1}
   s_{\overline{a_{1}^{k}a_{n}^{\ell}a_{i}a_{j}a_{i}^{-1}a_{j}^{-1}},a_{j}}^{-1}.
 \end{equation}
 where $1 \le i < j \le n$, $|i-j|\ge 2$.  Taking $i=1$ and $3\le j
 \le n-1$ gives: for $\ell=0$ (using (\ref{B-relation-type-1}));
 \begin{equation}
   s_{a_{1}^{k+1},a_j} =s_{a_{1}^{k},a_j},
 \end{equation}
 so by induction on $k$,
 \begin{equation}
 \label{B-relation-type-2-eliminate-k}
   s_{a_{1}^{k},a_j}=s_{1,a_j} \text{~~~for $3\le j \le n-1$},
 \end{equation}
 and for $\ell\ne0$;
 \begin{equation}
 \label{B-relation-type-2-1}
   s_{a_{1}^{k}a_{n}^{\ell},a_1}s_{a_{1}^{k+1}a_{n}^{\ell},a_j}
   s_{a_{1}^{k+1}a_{n}^{\ell},a_1}^{-1}s_{a_{1}^{k}a_{n}^{\ell},a_j}^{-1}.
 \end{equation}
 We will come back to relation (\ref{B-relation-type-2-1}) in a
 bit.

 Taking $i=1$ and $j=n$ in (\ref{B-relation-type-2})
 (and using (\ref{B-relation-type-1})) gives
 \begin{equation}
 \label{B-relation-type-2-2}
    s_{a_{1}^{k}a_{n}^{\ell},a_1}s_{a_{1}^{k}a_{n}^{\ell+1},a_1}^{-1}.
 \end{equation}
 So, by induction on $\ell$ (and (\ref{B-relation-type-1})) we get
 \begin{equation}
 \label{B-relation-trivial-for-1}
   s_{a_{1}^{k}a_{n}^{\ell},a_1}=1 \text{~~~for $k,\ell \in \Z$}.
 \end{equation}

 Taking $2\le i \le n-2$, $i+2\le j \le n$ in
 (\ref{B-relation-type-2}) gives
 \begin{eqnarray}
 \label{B-relation-type-2-3}
  \begin{cases}
    s_{a_{1}^{k}a_{n}^{\ell},a_i}s_{a_{1}^{k+1}a_{n}^{\ell},a_j}
    s_{a_{1}^{k+1}a_{n}^{\ell},a_i}^{-1}s_{a_{1}^{k}a_{n}^{\ell},a_j}^{-1}
    \text{~~~~~for $j\le n-1$}, \\
    s_{a_{1}^{k}a_{n}^{\ell},a_i}s_{a_{1}^{k}a_{n}^{\ell+1},a_i}^{-1}
    \text{~~~~~~~~~~~~~~~~~~~~~~~~~~~~~~~~~~~for $j=n$.}
  \end{cases}
 \end{eqnarray}
 In the case $j=n$ induction on $\ell$ gives
 \begin{equation}
   s_{a_{1}^{k}a_{n}^{\ell},a_{i}} =s_{a_{1}^{k},a_{i}} \text{~~~($2\le i
   \le n-2$)}.
 \end{equation}
 So from (\ref{B-relation-type-2-eliminate-k}) it follows
 \begin{eqnarray}
 \label{B-relation-type-2-3-result}
   s_{a_{1}^{k}a_{n}^{\ell},a_{i}} =
   \begin{cases}
     s_{1,a_{i}}  \text{~~~~~$3\le i \le n-2$} \\
     s_{a_{1}^{k},a_{2}} \text{~~~~~~~$i=2$}.
   \end{cases}
 \end{eqnarray}
 We come back to the case $j\le n-1$ later.

 Returning now to (\ref{B-relation-type-2-1}), we can use
 (\ref{B-relation-trivial-for-1}) to get
 \begin{equation*}
   s_{a_{1}^{k+1}a_{n}^{\ell},a_{j}}=s_{a_{1}^{k}a_{n}^{\ell}.a_{j}}
   \text{~~~~($3\le j \le n-1$)}.
 \end{equation*}
 Thus, by induction on $k$
 \begin{equation}
   s_{a_{1}^{k}a_{n}^{\ell},a_{j}} = s_{a_{n}^{\ell},a_{j}} \text{~~~~($3\le j \le n-1$)}.
 \end{equation}
 For $3\le j \le n-2$ we already know this (equation
 (\ref{B-relation-type-2-3-result})), so the only new information
 we get from (\ref{B-relation-type-2-1}) is
 \begin{equation}
   s_{a_{1}^{k}a_{n}^{\ell},a_{n-1}} = s_{a_{n}^{\ell},a_{n-1}}
   \text{~~~~($k\in \Z$)}.
 \end{equation}

 Collecting all the information we have obtained from
 $\tau(a_{1}^{k}a_{n}^{\ell}T_{i,j}(a_{1}^{k}a_{n}^{\ell})^{-1})$,
 $1\le i<j\le n,|i-j|\ge2$, we get:
 \begin{equation}
 \label{B-relation-type-2-final}
 \begin{split}
    s_{a_{1}^{k}a_{n}^{\ell},a_{1}} &= 1  \text{~~~($k,\ell\in \Z$)}, \\
    s_{a_{1}^{k}a_{n}^{\ell},a_{i}} &=
   \begin{cases}
     s_{1,a_{i}}  \text{~~~~~$3\le i \le n-2$}, \\
     s_{a_{1}^{k},a_{2}} \text{~~~~~~~$i=2$},
   \end{cases} \\
   s_{a_{1}^{k}a_{n}^{\ell},a_{n-1}} &= s_{a_{n}^{\ell},a_{n-1}},
 \end{split}
 \end{equation}
 and (from (\ref{B-relation-type-2-3})), for $2\le i \le n-3$ and $i+2\le j\le
 n-1$,
 \begin{equation}
  s_{a_{1}^{k}a_{n}^{\ell},a_i}s_{a_{1}^{k+1}a_{n}^{\ell},a_j}
  s_{a_{1}^{k+1}a_{n}^{\ell},a_i}^{-1}s_{a_{1}^{k}a_{n}^{\ell},a_j}^{-1}.
 \end{equation}
 This relation breaks up into the following cases (using
 (\ref{B-relation-type-2-final}))
 \begin{eqnarray}
 \label{B-relation-type-2-3-final}
   \begin{cases}
     s_{a_{1}^{k},a_{2}}s_{1,a_{j}}
     s_{a_{1}^{k+1},a_{2}}^{-1}s_{1,a_{j}}^{-1}
     \text{~~~~~~for $i=2$, $4\le j \le n-2$,} \\
     s_{a_{1}^{k},a_{2}}s_{a_{n}^{\ell},a_{n-1}}
     s_{a_{1}^{k+1},a_{2}}^{-1}s_{a_{n}^{\ell},a_{n-1}}^{-1}
     \text{~~~~~~for $i=2$, $j=n-1$,} \\
     s_{1,a_{i}}s_{1,a_{j}}
     s_{1,a_{i}}^{-1}s_{1,a_{j}}^{-1}
     \text{~~~~~~for $3\le i\le n-3$, $i+2\le j \le n-2$,} \\
     s_{1,a_{i}}s_{a_{n}^{\ell},a_{n-1}}
     s_{1,a_{i}}^{-1}s_{a_{n}^{\ell},a_{n-1}}^{-1}
     \text{~~~~~~for $3\le i\le n-3$, $j = n-1$,} \\
   \end{cases}
 \end{eqnarray}

 The third type of relations in (\ref{B-prelim-presentation});
 $\tau(a_{1}^{k}a_{n}^{\ell}R_{i}(a_{1}^{k}a_{n}^{\ell})^{-1})$,
 after rewriting using
 equation (\ref{reidemeister-rewriting-function}), are
 \begin{equation}
 \label{B-relation-type-3}
   s_{a_{1}^{k}a_{n}^{\ell},a_{i}}s_{a_{1}^{k+1}a_{n}^{\ell},a_{i+1}}
   s_{a_{1}^{k+2}a_{n}^{\ell},a_{i}}s_{a_{1}^{k+2}a_{n}^{\ell},a_{i+1}}^{-1}
   s_{a_{1}^{k+1}a_{n}^{\ell},a_{i}}^{-1}s_{a_{1}^{k}a_{n}^{\ell},a_{i+1}}^{-1},
 \end{equation}
 which break down as follows (using (\ref{B-relation-type-2-final})):
 \begin{eqnarray}
 \label{B-relation-type-3-final}
  \begin{cases}
    s_{a_{1}^{k+1},a_{2}}s_{a_{1}^{k+2},a_{2}}^{-1}
    s_{a_{1}^{k},a_{2}}^{-1}
    \text{~~~~~~ ($i=1$),} \\
    s_{a_{1}^{k},a_{2}}s_{1,a_{3}}
    s_{a_{1}^{k+2},a_{2}}s_{1,a_{3}}^{-1}
    s_{a_{1}^{k+1},a_{2}}^{-1}s_{1,a_{3}}^{-1}
    \text{~~~~~~($i=2$),} \\
    s_{1,a_{i}}s_{1,a_{i+1}}
    s_{1,a_{i}}s_{1,a_{i+1}}^{-1}
    s_{1,a_{i}}^{-1}s_{1,a_{i+1}}^{-1},
    \text{~~~~~~ for $3\le i \le n-3$,} \\
    s_{1,a_{n-2}}s_{a_{n}^{\ell},a_{n-1}}
    s_{1,a_{n-2}}s_{a_{n}^{\ell},a_{n-1}}^{-1}
    s_{1,a_{n-2}}^{-1}s_{a_{n}^{\ell},a_{n-1}}^{-1},
    \text{~~~~~~ ($i=n-2$),}
  \end{cases}
 \end{eqnarray}

 The fourth type of relations in (\ref{B-prelim-presentation});
 $\tau(a_{1}^{k}a_{n}^{\ell}R_{n-1}(a_{1}^{k}a_{n}^{\ell})^{-1})$, after rewriting using
 equation (\ref{reidemeister-rewriting-function}), is
 \begin{equation}
 \label{B-relation-type-4}
   s_{a_{n}^{\ell},a_{n-1}}s_{a_{n}^{\ell+1},a_{n-1}}
   s_{a_{n}^{\ell+2},a_{n-1}}^{-1}s_{a_{n}^{\ell+1},a_{n-1}}^{-1},
 \end{equation}
 where we have made extensive use of the relations
 (\ref{B-relation-type-2-final}).

 From (\ref{B-relation-type-2-final}) it follows that $\A_{B_n}'$ is
 generated by $s_{a_{1}^{k},a_{2}}$, $s_{1,a_{i}}$,
 and $s_{a_{n}^{\ell},a_{n-1}}$ for $k,\ell \in \Z$ and $3\le i \le
 n-2$.  For simplicity of notation let these generators be denoted
 by $p_k$, $q_i$, and $r_{\ell}$, respectively.  Thus, we have shown that the
 following is a set of defining relations for $\A_{B_n}'$:
 \begin{equation}
 \label{B-presentation}
 \begin{split}
   p_{k}q_{j}=q_{j}p_{k+1} &~ \text{~~~($4\le j \le n-2$, $k \in \Z$}), \\
   p_{k}r_{\ell}=r_{\ell}p_{k+1} &~ \text{~~~($k,\ell \in \Z$)}, \\
   q_{i}q_{j}=q_{j}q_{i} &~ \text{~~~($3\le i<j \le n-2$, $|i=j|\ge2$}), \\
   q_{i}r_{\ell}=r_{\ell}q_{i} &~ \text{~~~($3\le i \le n-3$)}, \\
   p_{k+1}p_{k+2}^{-1}p_{k}^{-1}  &~ \text{~~~($k\in \Z$)}, \\
   p_{k}q_{3}p_{k+2}q_{3}^{-1}p_{k+1}^{-1}q_{3}^{-1} &~  \text{~~~($k\in Z$)}, \\
   q_{i}q_{i+1}q_{i}=q_{i+1}q_{i}q_{i+1} &~ \text{~~~($3\le i \le n-3$)}, \\
   q_{n-2}r_{\ell}q_{n-2}=r_{\ell}q_{n-2}r_{\ell} &~ \text{~~~($\ell\in \Z$)}, \\
   r_{\ell}r_{\ell+1}r_{\ell+2}^{-1}r_{\ell+1}^{-1} &~ \text{~~~($\ell \in \Z$)},
 \end{split}
 \end{equation}
 The first four relations are from (\ref{B-relation-type-2-3-final}), the next four are
 from (\ref{B-relation-type-3-final}), and the last one is from (\ref{B-relation-type-4}).

 The fifth relation tells us that for $k \neq 0,1$, $p_k$ can be expressed in terms of $p_0$
 and $p_1$.  Similarly the last relation tells us that for $\ell \neq 0,1$, $r_{\ell}$ can be expressed
 in terms of $r_0$ and $r_1$.  From this it follows that $\A_{B_n}'$ is finitely
 generated.  Using lemmas \ref{lem:relation-replacement-1}
 and \ref{lem:relation-replacement-2} to replace the first, second
 and sixth relations, assuming we have added a new generator $b$
 and relation $b=p_0q_3p_{0}^{-1}$, we arrive at the following
 theorem.

 \begin{theorem}
  \label{thm:commutator-type-B-presentation}
  For every $n\geq3$ the commutator subgroup $\A_{B_n}'$ of the
  Artin group $\A_{B_n}$ is a finitely  generated group.
  Presentations for $\A_{B_n}'$, $n\ge 2$ are as follows: \\
  $\A_{B_2}'$ is a free group on countably many generators:
  \begin{gather*}
   [a_{2}^{\ell},a_{1}] ~~ (\ell\in \Z \setminus \{0,\pm1 \}),
   \quad [a_{1}^{k}a_{2},a_{1}] ~~ (k\in \Z \setminus \{0 \}).
  \end{gather*}
  $\A_{B_3}'$ is a free group on four generators:
  \begin{equation*}
   [a_{1}^{-1},a_{2}^{-1}] , \quad [a_{3},a_{2}][a_{1}^{-1},a_{2}^{-1}] ,
   \quad [a_{1},a_{2}][a_{1}^{-1},a_{2}^{-1}] , \quad
   [a_{1}a_{3},a_{2}][a_{1}^{-1},a_{2}^{-1}].
  \end{equation*}
  $\A_{B_4}'$ is the group generated by
  \begin{gather*}
    p_k=a_{1}^{k}a_{2}a_{1}^{-(k+1)}=[a_{1}^{k},a_{2}][a_{1}^{-1},a_{2}^{-1}],~~~
    (k\in \Z)\\
    q_{\ell}=a_{4}^{\ell}a_{3}(a_{1}a_{4}^{\ell})^{-1}=[a_{4}^{\ell},a_{3}][a_{2}^{-1},a_{3}^{-1}]
    [a_{1}^{-1},a_{2}^{-1}], ~~~(\ell\in \Z),
   \end{gather*}
   with defining relations
   \begin{gather*}
   p_{k+1}p_{k+2}^{-1}p_{k}^{-1}  \quad (k \in \Z), \\
   p_{k}q_{\ell}p_{k+2}=q_{\ell}p_{k+1}q_{\ell} \quad (k,\ell \in \Z), \\
   q_{\ell}q_{\ell+1}=q_{\ell+1}q_{\ell+2} \quad (3\le i \le n-3).
   \end{gather*}
  For $n\geq5$ the group $\A_{B_n}'$ is generated by
   \begin{gather*}
    p_0=a_{2}a_{1}^{-1},~~
    p_1=a_{1}a_{2}a_{1}^{-2},~~
    q_3=a_{3}a_{1}^{-1},~~
    r_{\ell}=a_{n}^{\ell}a_{n-1}(a_{1}a_{n}^{\ell})^{-1} ~~ (\ell \in \Z), \\
    b=a_{2}a_{1}^{-1}a_{3}a_{2}^{-1},~~
    q_i=a_{i}a_{1}^{-1} ~~(4\leq i \leq n-2),
   \end{gather*}
   with defining relations
   \begin{gather*}
   p_{0}q_{j}=q_{j}p_{1}, \quad  p_{1}q_{j}=q_{j}p_{o}^{-1}p_{1} \text{~~~($4\le j \le n-2$)}, \\
   p_{0}r_{\ell}=r_{\ell}p_{1}, \quad  p_{1}r_{\ell}=r_{\ell}p_{0}^{-1}p_{1} \text{~~~($\ell \in \Z$)}, \\
   q_{i}q_{j}=q_{j}q_{i}  \text{~~~($3\le i<j \le n-2$, $|i=j|\ge2$}), \\
   q_{i}r_{\ell}=r_{\ell}q_{i}  \text{~~~($3\le i \le n-3$)}, \\
   p_{0}q_{3}p_{0}^{-1}=b, \quad p_{0}bp_{0}^{-1}=b^2q_{3}^{-1}b, \\
   p_{1}q_3p_{1}^{-1} = q_{3}^{-1}b, \quad
   p_{1}bp_{1}^{-1}=(q_{3}^{-1}b)^{3}q_{3}^{-2}b, \\
   q_{i}q_{i+1}q_{i}=q_{i+1}q_{i}q_{i+1}  \text{~~~($3\le i \le n-3$)}, \\
   q_{n-2}r_{\ell}q_{n-2}=r_{\ell}q_{n-2}r_{\ell}  \text{~~~($\ell\in \Z$)}, \\
   r_{\ell}r_{\ell+1}r_{\ell+2}^{-1}r_{\ell+1}^{-1}  \text{~~~($\ell \in \Z$)},
   \end{gather*}
   \hfill $\Box$
 \end{theorem}

 \begin{corollary}
  \label{cor:commutator-type-B-perfect}
    For $n\ge5$ the commutator subgroup $\A_{B_n}'$ of the Artin
    group of type $B_n$ is finitely generated and perfect.
  \end{corollary}

  \begin{proof}
    Abelianizing the presentation of $\A_{B_n}'$ in the theorem
    results in a presentation of the trivial group.  Hence
    $\A_{B_n}''=\A_{B_n}'$.
  \end{proof}

\subsubsection{Type $D$}
\label{subsec:type-D}

The presentation of $\A_{D_n}$ is
 \begin{eqnarray*}
    \A_{D_n}=\la a_1,...,a_{n}: &~&a_{i}a_{j}=a_{j}a_{i}
    \text{~~~for $1\le i <j \le n-1$,$|i-j|\geq 2$},\\
    &~&a_{n}a_{j}=a_{j}a_{n} \text{~~~for $j \ne n-2$}, \\
    &~&a_{i}a_{i+1}a_{i}=a_{i+1}a_{i}a_{i+1}
    \text{~~~for $1 \leq i \leq n-2$ } \\
    &~&a_{n-2}a_{n}a_{n-2}=a_{n}a_{n-2}a_{n}
    \ra.
 \end{eqnarray*}

 As seen in figure \ref{fig:odd-coxeter-graphs} the graph
 $(D_n)_{odd}$ is connected. So by proposition~\ref{thm:connected-commutator}
 \begin{eqnarray*}
 \A_{D_n}'=\{U\in \A_{D_n} : \d(U)=0 \}.
 \end{eqnarray*}
 The computation of the presentation of $\A_{D_n}'$ is similar to
 that of $\A_{A_n}'$, so we will not include it.

 \begin{theorem}
 \label{thm:commutator-type-D-presentation}
   For every $n\geq4$ the commutator subgroup $\A_{D_n}'$ of the
   Artin group $\A_{D_n}$ is a finitely presented group.  $\A_{D_4}'$
   is the group generated by
   \begin{eqnarray*}
     p_0=a_{2}a_{1}^{-1},~~
     p_1=a_{1}a_{2}a_{1}^{-2},~~
     q_3=a_{3}a_{1}^{-1},~~\\
     q_4=a_{4}a_{1}^{-1},~~
     b=a_{2}a_{1}^{-1}a_{3}a_{2}^{-1},
     c=a_{2}a_{1}^{-1}a_{4}a_{2}^{-1},
   \end{eqnarray*}
   with defining relations
   \begin{gather*}
      b=p_0q_3p_{0}^{-1}, \quad
      p_0bp_{0}^{-1}=b^{2}q_3^{-1}b,\\
      p_1q_3p_{1}^{-1}=q_3^{-1}b, \quad
      p_1bp_{1}^{-1}=(q_3^{-1}b)^{3}q_3^{-2}b,\\
      c=p_0q_4p_{0}^{-1}, \quad
      p_0cp_{0}^{-1}=c^{2}q_4^{-1}c,\\
      p_1q_4p_{1}^{-1}=q_4^{-1}c, \quad
      p_1cp_{1}^{-1}=(q_{4}^{-1}c)^{3}q_{4}^{-2}c,\\
      q_{3}q_{4}=q_{4}q_{3}.
   \end{gather*}
   For $n\geq5$ the group $\A_{D_n}'$ is generated by
    \begin{eqnarray*}
       p_0=a_{2}a_{1}^{-1},\quad
       p_1=a_{1}a_{2}a_{1}^{-2},\quad\\
       q_{\ell}=a_{\ell}a_{1}^{-1}~~(3\leq \ell \leq n), \quad
       b=a_{2}a_{1}^{-1}a_{3}a_{2}^{-1},
    \end{eqnarray*}
   with defining relations
   \begin{gather*}
      b=p_0q_3p_{0}^{-1}, \quad
      p_0bp_{0}^{-1}=b^{2}q_3^{-1}b,\\
      p_1q_3p_{1}^{-1}=q_3^{-1}b, \quad
      p_1bp_{1}^{-1}=(q_3^{-1}b)^{3}q_3^{-2}b,\\
      p_{0}q_{j}=q_{j}p_{1}, \quad p_{1}q_{j}=q_{j}p_{0}^{-1}p_{1} \quad
      (4\le j \le n), \\
      q_iq_{i+1}q_i=q_{i+1}q_{i}q_{i+1}~~(3 \leq i\leq n-2), \\
      q_nq_{n-2}q_n=q_{n-2}q_{n}q_{n-2}, \\
      q_iq_j=q_jq_i \quad (3\le i <j \le n-1, |i-j|\ge2), \\
      q_nq_j=q_jq_n \quad (j\ne n-2).
   \end{gather*}
   \hfill $\Box$
  \end{theorem}

  \begin{corollary}
  \label{cor:commutator-type-D-perfect}
    For $n\ge5$ the commutator subgroup $\A_{D_n}'$ of the Artin
    group of type $D_n$ is finitely presented and perfect. \hfill $\Box$
  \end{corollary}

\subsubsection{Type $E$}

The presentation of $\A_{E_n}$, $n=6,7,$ or $8$, is
 \begin{eqnarray*}
    \A_{E_n}=\la a_1,...,a_{n}: &~&a_{i}a_{j}=a_{j}a_{i}
    \text{~~~for $1\le i <j \le n-1$,$|i-j|\geq 2$},\\
    &~&a_{i}a_{n}=a_{n}a_{i} \text{~~~for $i \ne 3$}, \\
    &~&a_{i}a_{i+1}a_{i}=a_{i+1}a_{i}a_{i+1}
    \text{~~~for $1 \le i \le n-2$ } \\
    &~&a_{3}a_{n}a_{3}=a_{n}a_{3}a_{n}
    \ra.
 \end{eqnarray*}

 As seen in figure \ref{fig:odd-coxeter-graphs} the graph
 $(E_n)_{odd}$ is connected. So by proposition~\ref{thm:connected-commutator}
 \begin{eqnarray*}
 \A_{E_n}'=\{U\in \A_{E_n} : \d(U)=0 \}.
 \end{eqnarray*}
 The computation of the presentation of $\A_{E_n}'$ is similar to
 that of $\A_{A_n}'$.

 \begin{theorem}
 \label{thm:commutator-type-E-presentation}
   For $n=6,7,$ or $8$ the commutator subgroup $\A_{E_n}'$ of the
   Artin group $\A_{E_n}$ is a finitely presented group.  $\A_{E_n}'$
   is the group generated by
   \begin{eqnarray*}
     p_0=a_{2}a_{1}^{-1}, \quad
     p_1=a_{1}a_{2}a_{1}^{-2}, \quad
     q_{\ell}=a_{\ell}a_{1}^{-1} ~~(3\le \ell \le n), \quad
     b=a_{2}a_{1}^{-1}a_{3}a_{2}^{-1},
   \end{eqnarray*}
   with defining relations
   \begin{gather*}
      b=p_0q_3p_{0}^{-1}, \quad
      p_0bp_{0}^{-1}=b^{2}q_3^{-1}b,\\
      p_1q_3p_{1}^{-1}=q_3^{-1}b, \quad
      p_1bp_{1}^{-1}=(q_3^{-1}b)^{3}q_3^{-2}b,\\
      p_{0}q_{j}=q_{j}p_{1}, \quad p_1q_j=q_{j}p_{0}^{-1}p_{1} \quad (4\le
      j \le n), \\
      q_{i}q_{i+1}q_{i}=q_{i+1}q_{i}q_{i+1} \quad (3\le i \le n-2), \\
      q_{n}q_{3}q_{n}=q_{3}q_{n}q_{3},  \\
      q_{i}q_{j}=q_{j}q_{i} \quad (3\le i <j \le n-1, |i-j|\ge 2), \\
      q_{i}q_{n}=q_{n}q_{i}\quad (4\le i \le n-1).
   \end{gather*}
   \hfill $\Box$
  \end{theorem}

  \begin{corollary}
  \label{cor:commutator-type-E-perfect}
    For $n=6,7,$ or $8$ the commutator subgroup $\A_{E_n}'$ of the Artin
    group of type $E_n$ is finitely presented and perfect. \hfill $\Box$
  \end{corollary}

\subsubsection{Type $F$}
\label{subsec:type-F}

The presentation of $\A_{F_4}$ is
 \begin{eqnarray*}
    \A_{F_4}=\la a_1,a_2,a_3,a_{4}: &~&a_{i}a_{j}=a_{j}a_{i}
    \text{~~~for $|i-j|\geq 2$},\\
    &~&a_{1}a_{2}a_{1}=a_{2}a_{1}a_{2}, \\
    &~&a_{2}a_{3}a_{2}a_{3}=a_{3}a_{2}a_{3}a_{2}, \\
    &~&a_{3}a_{4}a_{3}=a_{4}a_{3}a_{4}
    \ra.
 \end{eqnarray*}

 As seen in figure \ref{fig:odd-coxeter-graphs} the graph
 $(E_n)_{odd}$ has two components: $\Gam_1$ and
 $\Gam_2$, where $\Gam_1$ denotes the component containing the
 vertices $a_1,a_2$, and $\Gam_2$ the component containing the
 vertices $a_3,a_4$.  Let $\d_1$ and $\d_2$ denote the
 associated degree maps, respectively, so from proposition
 ~\ref{thm:characterization-commuator-subgroup}
 \begin{eqnarray*}
 \A_{F_4}'=\{U\in \A_{F_4} : \d_1(U)=0 \text{~and~} \d_2(U)=0
 \}.
 \end{eqnarray*}

 By a computation similar to that of $B_n$ we get the following.
 \begin{theorem}
  \label{thm:commutator-type-F-presentation}
   The commutator subgroup $\A_{F_4}'$ of the Artin group of type
   $F_4$ is the group generated by
   \begin{gather*}
     p_k=a_{1}^{k}a_{2}a_{1}^{-(k+1)}=[a_{1}^{k},a_{2}][a_{1}^{-1},a_{2}^{-1}]
     \quad (k\in \Z), \\
     q_{\ell}=a_{4}^{\ell}a_{3}a_{4}^{-(\ell+1)}=[a_{4}^{\ell},a_{3}][a_{4}^{-1},a_{3}^{-1}]
     \quad (\ell\in \Z),
   \end{gather*}
   with defining relations
   \begin{gather*}
    p_{k} =  p_{k+1}p_{k+2}^{-1}~~ (k\in \Z), \quad q_{\ell} = q_{\ell+1}q_{\ell+2}^{-1}~~(\ell\in
     \Z),\\
     p_{k}q_{\ell}p_{k+1}q_{\ell+1}=q_{\ell}p_{k}q_{\ell+1}p_{k+1}~~(k,\ell \in \Z).
   \end{gather*}
 \end{theorem}

 The first two types of relations in the above presentation  tell
 us that for $k\ne 0,1$, $p_k$ can be expressed in terms of $p_0$
 and $p_1$, and similarly for $q_{\ell}$.  Thus $\A_{F_4}'$ is finitely
 generated.  However, $\A_{F_4}'$ is not perfect since abelianizing
 the above presentation gives $\A_{F_4}'/\A_{F_4}'' \simeq \Z^4$.

\subsubsection{Type $H$}
\label{subsec:type-H}

The presentation of $\A_{H_n}$, $n=3$ or $4$, is
 \begin{eqnarray*}
    \A_{H_n}=\la a_1,...,a_{n}: &~&a_{i}a_{j}=a_{j}a_{i}
    \text{~~~for $|i-j|\geq 2$},\\
    &~&a_{1}a_{2}a_{1}a_{2}a_{1}=a_{2}a_{1}a_{2}a_{1}a_{2}, \\
    &~&a_{i}a_{i+1}a_{i}=a_{i+1}a_{i}a_{i+1}
    \text{~~~for $2 \le i \le n-1$ }
    \ra.
 \end{eqnarray*}

 As seen in figure \ref{fig:odd-coxeter-graphs} the graph
 $(H_n)_{odd}$ is connected. So by proposition~\ref{thm:connected-commutator}
 \begin{eqnarray*}
 \A_{H_n}'=\{U\in \A_{H_n} : \d(U)=0 \}.
 \end{eqnarray*}
 The computation of the presentation of $\A_{H_n}'$ is similar to
 that of $\A_{A_n}'$.

 \begin{theorem}
 \label{thm:commutator-type-H-presentation}
   For $n=3$ or $4$ the commutator subgroup $\A_{H_n}'$ of the
   Artin group $\A_{H_n}$ is the group generated by
   \begin{eqnarray*}
     p_k=a_{1}^{k}a_{2}a_{1}^{-(k+1)} \quad (k\in \Z), \quad
     q_{\ell}=a_{\ell}a_{1}^{-\ell} \quad (3\le \ell \le n),
   \end{eqnarray*}
   with defining relations
   \begin{gather*}
      p_{k}q_{j}=q_{j}p_{k+1} ~~(4\le j \le n), \\
      p_{k+1}p_{k+3}p_{k+4}^{-1}p_{k+2}^{-1}p_{k}^{-1} ~~ (k \in \Z), \\
      p_{k}q_{3}p_{k+2}q_{3}^{-1}p_{k+1}^{-1}q_{3}^{-1} \\
      q_{i}q_{i+1}q_{i}=q_{i+1}q_{i}q_{i+1} \quad (3\le i \le n-1).
   \end{gather*}
   \hfill $\Box$
  \end{theorem}

   The second relation tells us that for $k\ne 0,1,2,3$, $p_k$ can
   be expressed in terms of $p_0,p_1,p_2,$ and $p_3$.  Thus,
   $\A_{H_n}'$ is finitely generated.  Abelianizing the above
   presentation results in the trivial group.  Thus, we have the
   following.

  \begin{corollary}
  \label{cor:commutator-type-H-perfect}
    For $n=3$ or $4$ the commutator subgroup $\A_{H_n}'$ of the Artin
    group of type $H_n$ is finitely generated and perfect. \hfill $\Box$
  \end{corollary}

\subsubsection{Type $I$}
\label{subsec:type-I}

 The presentation of $I_2(m)$, $m\ge 5$, is
 \begin{eqnarray*}
    \A_{I_2(m)}=\la a_1,a_2: \la a_{1}a_{2}\ra^m= \la
    a_{2}a_{1}\ra^m \ra.
 \end{eqnarray*}

 In figure \ref{fig:odd-coxeter-graphs} the graph
 $(I_2(m))_{odd}$ is connected for $m$ odd and disconnected for
 $m$ even.  Thus, different computations must be done for these
 two cases.  We have the following.

 \begin{theorem}
 \label{thm:commutator-type-I-presentation}
   The commutator subgroup $\A_{I_2(m)}'$ of the Artin group of
   type $I_{2}(m)$, $m\ge 5$, is the free group generated by the $(m-1)$-generators
   \begin{gather*}
     a_{1}^{k}a_{2}a_{1}^{-(k+1)} \quad (k\in\{0,1,2,\ldots, m-2
     \}),
   \end{gather*}
   when $m$ is odd, and is the free group with countably many
   generators
   \begin{gather*}
     [a_{2}^{\ell},a_{1}] ~~(\ell\in \Z \setminus \{-(m/2-1) \}),
     \quad [a_{1}^{j}a_{2}^{\ell},a_{1}] ~~(\ell\in \Z,~~j=1,2,\ldots,m/2-3) , \\
     [a_{1}^{m/2-2}a_{2}^{\ell},a_{1}] ~~(\ell\in \Z \setminus \{m/2-1 \}),
     \quad [a_{1}^{k}a_{2},a_{1}] ~~(k \in \Z).
   \end{gather*}
   when $m$ is even.
 \end{theorem}

\subsubsection{Summary of results}
\label{subsec:summary}

 Table \ref{tab:artin-commutator-subgroups} summarizes the results in this section.
 The question marks (?) in the table indicate that it is unknown
 whether the commutator subgroup is finitely presented.  However, we do
 know that for these cases the group is finitely generated.  If
 one finds more general relation equivalences along the
 lines of lemmas \ref{lem:relation-replacement-1} and \ref{lem:relation-replacement-2}
 then we may be able to show that these groups are indeed finitely
 presented.

 \begin{table}[top]
  \begin{center}
    \begin{tabular}{|c|c|c|}
    \hline
     Type $\Gamma$ & finitely generated/presented & perfect                   \\ \hline \hline
         $A_n$       &  yes/yes & $n=1,2,3$ : no,             \\
                     &          & $n\ge 4$ : yes \\ \hline
                     & $n=2$ : no, $n\ge 3$ : yes & $n=2,3,4$ : no, \\
         $B_n$       & / & $n\ge5$ : yes \\
                     &  $n=3$ : yes, $n\ge 3$ : ? & \\ \hline
         $D_n$       &  yes/yes & $n=4$ : no,     \\
                     &          & $n\ge 5$ : yes \\ \hline
         $E_n$       &  yes/yes & yes        \\ \hline
         $F_4$       &  yes/? & no         \\ \hline
         $H_n$       &  yes/? & yes        \\ \hline
         $I_{2}(m)$ ($m$ even) & no/no& no    \\
         ~~~~~~~~~($m$ odd)   & yes/yes  & no          \\ \hline
    \end{tabular}
  \end{center}
  \caption{Properties of the commutator subgroups}
  \label{tab:artin-commutator-subgroups}
 \end{table}

\section{Local indicability}
\subsection{Definitions and generalities}
\label{sec4.1}

 A group $G$ is {\bf indicable} if there exists a {\it nontrivial}
 homomorphism $G\longrightarrow \Z$ (called an indexing
 function). A group $G$ is {\bf locally indicable}\index{locally indicable} if every
 nontrivial, finitely generated subgroup is indicable.  Notice, finite
 groups cannot be indicable, so locally indicable groups must be torsion-free.

 Local indicability was introduced by Higman \cite{Hi40b} in connection with the zero-divisor
 conjecture, which asserts that if $G$ is a torsion-free group, then its integral group ring
 $\Z G$ has no zero divisors.  Although still unsolved, the conjecture is true for groups
 which are locally indicable.

 Every free group is locally indicable.  Indeed, it is well known
 that every subgroup of a free group is itself free, and since free
 groups are clearly indicable the result follows.

 Local indicability is clearly inherited by subgroups.
 The following simple theorem shows that the category of locally
 indicable groups is preseved under extensions.
 \begin{theorem}
 \label{thm:local-indicability-extensions}
 If $K,H$ and $G$ are groups such that $K$ and $H$ are locally
 indicable and fit into a short exact sequence
 $$
 \begin{CD}
   1 \longrightarrow K @>{\phi}>> G  @>{\psi}>> H \longrightarrow 1,
 \end{CD}
 $$
 then G is locally indicable.
 \end{theorem}

 \begin{proof}
   Let $g_1,\ldots,g_n \in G$, and let $\la g_1,\ldots,g_n \ra$
   denote the subgroup of $G$ which they generate.  If
   $\psi(\la g_1,\ldots,g_n \ra)\ne \{1 \}$ then by the local
   indicability of $H$ there exists a nontrivial homomorphism
   $f:\psi(\la g_1,\ldots,g_n \ra) \longrightarrow \Z$.  Thus, the
   map
   \begin{equation*}
   f\circ\psi:\la g_1,\ldots,g_n \ra \longrightarrow \Z
   \end{equation*}
   is
   nontrivial.  Else, if $\psi(\la g_1,\ldots,g_n \ra) = \{1 \}$
   then $g_1,\ldots,g_n \in \ker \psi = \im \phi$ (by exactness),
   so there exist $k_1,\ldots,k_n \in K$ such that
   $\phi(k_i)=g_i$, for all $i$.  Since $\phi$ is
   one-to-one (short exact sequence) then
   $\phi:\la k_1,\ldots,k_n \ra\longrightarrow\la g_1,\ldots,g_n
   \ra$
   is an isomorphism.  By the local indicability of
   $K$ there exists a nontrivial homomorphism
   $h:\la k_1,\ldots,k_n \ra \longrightarrow \Z$, therefore the
   map
   \begin{equation*}
     h\circ\phi^{-1}: \la g_1,\ldots,g_n \ra \longrightarrow \Z
   \end{equation*}
   is nontrivial.
 \end{proof}

  \begin{corollary}
  \label{cor:direct-product-li}
    If $G$ and $H$ are locally indicable then so is $G\oplus H$.
  \end{corollary}

  \begin{proof}
    The sequence
    $$
    \begin{CD}
     1 \longrightarrow H @>{\phi}>> G\oplus H @>{\psi}>> G
     \longrightarrow 1
    \end{CD}
    $$
    where $\phi(h)=(1,h)$ and $\psi(g,h)=g$
    is exact, so the theorem applies.
  \end{proof}

  If $G$ and $H$ are groups and $\phi : G \longrightarrow {\rm
  Aut}(H)$.  The {\bf semidirect product}\index{semidirect product} of $G$ and $H$ is
  defined to be the set $H \times G$ with binary operation
  \begin{equation*}
   (h_1,g_1)\cdot(h_2,g_2) =(h_1\cdot g_1\ast h_2,g_1g_2)
  \end{equation*}
  where $g\ast h$ denotes the action of $G$ on $H$ determined by $\phi$,
  i.e. $g\ast h:= \phi(g)(h)\in H$.  This group is denoted by $H \rtimes_{\phi}
  G$.

  \begin{corollary}
  \label{cor:semi-direct-product-li}
    If $G$ and $H$ are locally indicable then so is $H \rtimes_{\phi}
    G$.
  \end{corollary}

  \begin{proof}
    If $\psi: H\rtimes_{\phi} G \longrightarrow G$ denotes the map
    $(h,g)\longmapsto g$ then $\ker \psi = H$ and the groups fit
    into the exact sequence
    $$
    \begin{CD}
     1 \longrightarrow H @>{{\rm incl.}}>> H \rtimes_{\phi}
     G  @>{\psi}>> G
     \longrightarrow 1
    \end{CD}
    $$
  \end{proof}

  The following theorem of Brodskii \cite{Br80}, \cite{Br84},
  which was discovered independently by Howie \cite{Ho82}, \cite{Ho00},
  tells us that the class of torsion-free $1$-relator groups lies
  inside the class of locally indicable groups.  Also, for
  $1$-relator groups: locally indicable $\Leftrightarrow$ torsion
  free.
  \begin{theorem}
  \label{thm:Howie}
   A torsion-free $1$-relator group is locally indicable.
  \end{theorem}

  To show a group is not locally indicable we need to show there
  exists a finitely generated subgroup in which the only
  homomorphism into $\Z$ is the trivial homomorphism.

  \begin{theorem}
  \label{thm:fg-perfect-sg-not-li}
    If $G$ contains a finitely generated perfect sugroup then $G$
    is not locally indicable.
  \end{theorem}
  \begin{proof}
    The image of a commutator $[a,b]:=aba^{-1}b^{-1}$ under
    a homomorphism into $\Z$ is 0, thus the image of a perfect
    group is trivial.
  \end{proof}

\subsection{Local indicability of spherical Artin groups}
\label{sec:local-indicability}

 Since spherical-type Artin groups are torsion-free,
 theorem \ref{thm:Howie}
 implies that the Artin groups of type $A_2,B_2$, and $I_2(m)$
 ($m\ge5$) are locally indicable.  In this section we determine the local
 indicability of all\footnote{with the exception of type $F_4$
 which at this time remains undetermined.}
 irreducible spherical-type Artin groups.

 It is of interest to note that the discussion in section
 \ref{sec:characterization}, in particular proposition
 \ref{thm:characterization-commuator-subgroup}, shows that
 an Artin group $\Ar$ and its commutator subgroup $\Ar'$ fit
 into a short exact sequence:
 $$
 \begin{CD}
  1 \longrightarrow \Ar' \longrightarrow \Ar @>{\d_{\Gam}}>> \Z^m
  \longrightarrow 1,
 \end{CD}
 $$
 where $m$ is the number of connected components in $\Gam_{odd}$,
 and $\d_{\Gam}$ is the degree map defined in \ref{sec:characterization}, which can be
 identified with the abelianization map.  Thus,
 the local indicability of an Artin group $\Ar$ is completely determined
 by the local indicability of its commutator subgroup $\Ar'$ (by theorem
 \ref{thm:local-indicability-extensions}).  In other words,
 \begin{eqnarray*}
  \Ar \text{ is locally indicable} \Longleftrightarrow \Ar' \text{ is locally indicable.}
 \end{eqnarray*}
  This gives another proof that
 the Artin groups of type $A_2,B_2$, and $I_2(m)$
 ($m\ge5$) are locally indicable, since their corresponding
 commutator subgroups are free groups, as already shown.
\subsubsection{Type $A$}
\label{subsec:type-A-indicability}

$A_{A_1}$ is clearly locally indicable since $A_{A_1}\simeq \Z$,
and, as noted above, $\A_{A_2}$ is also locally indicable.

 For $\A_{A_3}$, theorem \ref{thm:commutator-subgroup-A3}
 tells us $\A_{A_3}'$ is the semidirect
 product of two free groups, thus $\A_{A_3}'$ is locally indicable.
 It follows from our remarks above that $\A_{A_3}$ is also locally
 indicable.


 As for $\A_{A_n}$, $n\ge 4$, corollary
 \ref{cor:commutator-type-A-perfect} and theorem
 \ref{thm:fg-perfect-sg-not-li} imply that $\A_{A_n}$ is not
 locally indicable.

 Thus, we have the following theorem.
 \begin{theorem}
 \label{thm:type-A-local-indicability}
  $\A_{A_n}$ is locally indicable if and only if $n=1,2,$ or $3$.
 \end{theorem}

\subsubsection{Type $B$}
\label{subsec:type-B-indicability}

We saw above $\A_{B_2}$ is locally indicable.    For $n=3$ and $4$
we argue as follows.

Let $\P_{n+1}^{n+1}$ denote the $(n+1)$-pure braids in $\Br_{n+1}
= \A_{A_n}$, that is the braids which only permute the first
$n$-strings.  Letting $b_1,\ldots,b_n$ denote the generators of
$\A_{B_n}$ a theorem of Crisp \cite{Cr99} states
\begin{theorem}
  The map
  \begin{equation*}
    \phi : \A_{B_n} \longrightarrow \A_{A_n}
  \end{equation*}
  defined by
  \begin{equation*}
    b_i \longmapsto a_i, \quad b_n \longmapsto a_{n}^{2}
  \end{equation*}
  is an injective homomorphism onto $\P_{n+1}^{n+1}$.  That is,
  $\A_{B_n}\simeq \P_{n+1}^{n+1} < \Br_{n+1} = \A_{A_n}$.
\end{theorem}

By "forgetting the $n^{\rm th}$-strand" we get a homomorphism
$f:\P_{n+1}^{n+1} \longrightarrow \Br_n$ which fits into the short
exact sequence $$
\begin{CD}
  1 \longrightarrow K \longrightarrow \P_{n+1}^{n+1} @>{f}>>
  \Br_n \longrightarrow 1,
\end{CD}
$$ where $K=\ker~f= \{\beta \in \P_{n+1}^{n+1}: $ the first $ n $
strings of $ \beta $ are trivial$ \}$.  It is known that $K \simeq
F_n$, the free group of rank $n$.  Since $F_n$ is locally
indicable and $\Br_n$ $(n=3,4)$ is locally indicable then so is
$\A_{B_n}$, for $n=3,4$.  Futhermore, the above exact sequence is
actually a split exact sequence so $\A_{B_n}\simeq \P_{n+1}^{n+1}
\simeq F_n \rtimes \Br_n$.

As for $\A_{B_n}$, $n\ge 5$, corollary
\ref{cor:commutator-type-B-perfect} and theorem
\ref{thm:fg-perfect-sg-not-li} imply that $\A_{B_n}$ is not
locally indicable, for $n\ge 5$.

Thus, we have the following theorem.
 \begin{theorem}
  $\A_{B_n}$ is locally indicable if and only if $n \le 4$.
 \end{theorem}
\subsubsection{Type $D$}
\label{subsec:type-D-indicability}

It follows corollary \ref{cor:commutator-type-D-perfect} and
\ref{thm:fg-perfect-sg-not-li} that $\A_{D_n}$ is not locally
indicable for $n \ge 5$.  As for $\A_{D_4}$, we will show it is
locally indicable as follows.

A theorem of Crisp and Paris \cite{CP02} says:

\begin{theorem}
\label{thm:D-semidirect-of-free-and-A}
  Let $F_{n-1}$ denote the free group of rank $n-1$.  There is an
  action $\rho:\A_{A_{n-1}} \longrightarrow {\Aut}(F_{n-1})$ such
  that $\A_{D_n} \simeq F_{n-1} \rtimes \A_{A_{n-1}}$ and $\rho$
  is faithful.
\end{theorem}

Since $\A_{A_3}$ and $F_{3}$ are locally indicable, then so is
$\A_{D_4}$.  Thus, we have the following theorem.
 \begin{theorem}
  $\A_{D_n}$ is locally indicable if and only if $n=4$.
 \end{theorem}

\subsubsection{Type $E$}
\label{subsec:type-E-indicability}

Since the commutator subgroups of $\A_{E_n}$, $n=6,7,8$, are
finitely generated and perfect (corollary
\ref{cor:commutator-type-E-perfect}) then $\A_{E_n}$ is not
locally indicable.

\subsubsection{Type $F$}
\label{subsec:type-F-indicability}

Unfortunately, we have yet to determine the local indicability of
the Artin group $\A_{F_4}$.


\subsubsection{Type $H$}
\label{subsec:type-H-indicability}

Since the commutator subgroups of $\A_{H_n}$, $n=3,4$, are
finitely generated and perfect (corollary
\ref{cor:commutator-type-H-perfect}) then $\A_{H_n}$ is not
locally indicable.

\subsubsection{Type $I$}
\label{subsec:type-I-indicability}


 As noted above, since the commutator subgroup $\A_{I_2(m)}'$
 of $\A_{I_2(m)}$
 $(m\ge 5)$ is a free group
 (theorem~\ref{thm:commutator-type-I-presentation}) then
 $\A_{I_2(m)}'$ is locally indicable and therefore so
 is $\A_{I_2(m)}$.  One could also apply theorem
 \ref{thm:Howie} to conclude the same result.


\section{Open questions: Orderability}

In this section we discuss the connection between the theory of
orderable groups and the theory of locally indicable groups.
Then we discuss the current state of the orderability of the
irreducible spherical-type Artin groups.

\subsection{Orderable Groups}
\label{sec:orderable-groups}

A group or monoid $G$ is {\bf right-orderable}
\index{orderable!right} if there exists a strict linear ordering
$<$ of its elements which is right-invariant: $g<h$ implies
$gk<hk$ for all $g,h,k$ in $G$.  If there is an ordering of $G$
which is invariant under multiplication on both sides, we say that
$G$ is {\bf orderable} \index{orderable} or for emphasis {\bf
bi-orderable} \index{orderable!bi-orderable}.

\begin{theorem}
  $G$ is right-orderable if and only if there exists a subset $\P
  \subset G$ such that:
  \begin{center}
    $\P\cdot \P \subset \P$  (subsemigroup), \\
    $G \setminus \{ 1 \} = \P \sqcup \P^{-1}$.
  \end{center}
\end{theorem}
\begin{proof}
  Given $\P$ define $<$ by: $g < h$ iff $hg^{-1}\in \P$.
  Given $<$ take $\P=\{g\in G: 1<g \}$.
\end{proof}

The ordering is a bi-ordering if and only if $\P$ exists as above and also
\begin{center}
 $g\P g^{-1} \subset \P, \quad \forall g \in G$.
\end{center}
 The set $\P \subset G$ in the previous theorem is called the
{\bf positive cone} \index{positive cone} with respect to the
ordering $<$.

The theory of orderable groups is well over a hundred years old.
For a general exposition on the theory of orderable groups see
\cite{MR77} or \cite{KK74}.  We will list just a few properties of
orderable groups.

A group is right-orderable if and only if it is left-orderable (by a possibly
different ordering.
The class of right-orderable groups is closed under: subgroups,
direct  products, free products, semidirect products,
and extension.  The class of orderable groups is closed under:
subgroups, direct products, free products, but not necessarily
extensions.  Both right-orderability and bi-orderability are local
properties: a group has the property if and only if every
finitely-generated subgroup has it.

Knowing a group is right-orderable or bi-orderable provides useful
information about the internal structure of the group.  For
example, if $G$ is right-orderable then it must be torsion-free:
for $1<g$ implies $g<g^2<g^3<\cdots<g^n<\cdots$.  Moreover, if $G$
is bi-orderable then $G$ has no {\bf generalised torsion}
(products of conjugates of a nontrivial element being trivial),
$G$ has unique roots: $g^n=h^n \Rightarrow g=h$, and if
$[g^n,h]=1$ in $G$ then $[g,h]=1$.  Further consequences of
orderablility are as follows. For any group $G$, let $\Z G$ denote
the {\bf group ring} \index{group ring} of formal linear
combinations $n_1g_1+\cdots n_kg_k$.
\begin{theorem}
\label{thm:zdc}
  If $G$ is right-orderable, then $\Z G$ has no zero divisors. \qed
\end{theorem}
\begin{theorem}
 (Malcev, Neumann) If $G$ is bi-orderable, then $\Z G$ embeds in a
 division ring. \qed
\end{theorem}
\begin{theorem}
 (LaGrange, Rhemtulla)  If $G$ is right orderable and $H$ is any
 group, then $\Z G\simeq \Z H$ implies $G\simeq H$. \qed

\end{theorem}

 The following theorems give the connection between orderable groups
and locally indicable groups (see \cite{RR02}).  The first is due
to Levi \cite{Lev}, and the second due to Burns and Hale \cite{BH} .

 \begin{theorem}
 \label{thm:biord-imp-li-imp-ro}
Every bi-orderable group is locally indicable. \qed
  \end{theorem}

The Artin group of type $A_2$, for example, shows the converse does not hold.

\begin{theorem}
\label{liro}
 Every locally indicable group is right-orderable. \qed
  \end{theorem}

Bergman \cite{Ber} was the first to publish examples demonstrating
that the converse to theorem \ref{liro} is false. Note that the
Artin groups of type $A_n,\/ n \ge 4$ are also examples of
right-orderable groups which are not locally indicable.

 One final connection between local indicability and right-orderability
 was given by Rhemtulla and Rolfsen \cite{RR02}.

 \begin{theorem}
 \label{thm:Rhemtulla-Rolfsen}
   (Rhemtulla, Rolfsen)
   Suppose $(G,<)$ is right-ordered and there is a finite-index
   subgroup $H$ of $G$ such that $(H,<)$ is a bi-ordered group.
   Then $G$ is locally indicable.
 \end{theorem}

 An application of this theorem is as follows.  It is known that the braid groups $\Br_n=\A_{A_{n-1}}$ are right orderable
 \cite{DDRW02} and that the pure braids $\P_n$ are bi-orderable \cite{KRo02}.
 However, theorem \ref{thm:type-A-local-indicability} tells us
 that $\Br_n$ is not locally indicable for $n\ge 5$ therefore,
 by theorem \ref{thm:Rhemtulla-Rolfsen}, the
 bi-ordering on $\P_n$ and the right-ordering
 on $\Br_n$ are incompatible for $n\ge5$.

\subsection{Ordering spherical Artin groups}
\label{sec:orderability-spherical-type-artin}

The first proof the that braid groups $\Br_n$ enjoy a
right-invariant total ordering was given in \cite{De92},
\cite{De94}. Since then several quite different approaches have
been applied to understand this phenomenon.\footnote{For a
comprehensive look at this problem and all the different
approaches used to understand it see the book \cite{DDRW02}.} The
following theorems summarize the state of our knowledge regarding
orderability of the spherical Artin groups.

\begin{theorem}
\label{artino}
The Artin groups of type $A_n (n \ge 2), B, D, E, F, H$ and $I$ are not bi-orderable.
\end{theorem}

\begin{proof}
The Artin group of type $A_2$ is not biorderable, as it does not have unique roots: $aba=bab$ implies that $(ab)^3 = (ba)^3$, whereas $ab \ne ba$ (their images in the Coxeter group are distinct).  For similar reasons, the Artin groups of type $B_2$ and $I_2(m)$ are not biorderable.  The others listed all contain a type $A_2$ Artin group as a (parabolic) subgroup, and therefore cannot be bi-ordered.
\end{proof}

In other words, only the simplest irreducible spherical Artin group $\A_{A_1} \cong \Z$ can be bi-ordered.  On the other hand, many of them (perhaps all those of spherical type) have a right-invariant ordering.

\begin{theorem}
\label{artinro}
The Artin groups of type $A, B, D$ and $I$ are right-orderable.
\end{theorem}

\begin{proof}
Type $A$ are the braid groups, shown to be right-orderable by Dehornoy.  Since the Artin groups of type $B$ and $I$ embed in type $A$, they are also right-orderable (alternatively, type $I$ are right-orderable because they are locally indicable).  Finally, theorem \ref{thm:D-semidirect-of-free-and-A} and the fact that right-orderability is closed under extensions, implies the right-orderability of type $D$ Artin groups.
\end{proof}

However, we do not know whether
the remaining irreducible spherical-type Artin groups are right-orderable.
One approach is to reduce the problem to showing that the positive
Artin monoid is right-orderable.  If $\Gam$ is a Coxeter graph, we define the
corresponding Artin monoid $\pAr$ to be the monoid with the same presentation as
the Artin group, but interpreted as a monoid presentation.  That is,
\begin{eqnarray*}
  \pAr = \langle a \in \Sigma : \langle ab \rangle^{m_{ab}}=\langle ba
  \rangle^{m_{ab}} \text{ if } m_{ab}< \infty \rangle^+.
\end{eqnarray*}

It is known \cite{Pa01} that $ \pAr$ injects in $\Ar$ as the submonoid of words in the canonical
generators with no negative exponents.

\subsubsection{Ordering the Monoid is Sufficient}
\label{subsec:orderable-monoid}

We will show that for a Coxeter graph $\Gam$ of spherical type  the
Artin group $\Ar$ is right-orderable (resp. bi-orderable) if and
only if the Artin monoid $\pAr$ is right-orderable (resp.
bi-orderable). One direction is of course trivial.

Let $\Ar$ be an Artin group of spherical-type.  Brieskorn and
Saito \cite{BS72}, generalizing a result of Garside, have shown:

{\it For each $U\in \Ar$ there exist $U_1,U_2\in \pAr$, where
$U_2$ is central in $\Ar$, such that}
\begin{center}
  $U=U_1U_{2}^{-1}$.
\end{center}
All decompositions of elements of $\Ar$ in this section are
assumed to be of this form.

Suppose $\pAr$ is right-orderable, let $<^+$ be such a
right-invariant linear ordering.  We wish to prove that $\Ar$ is
right-orderable.

The following lemma indicates how we should extend the ordering on
the monoid to the entire group.
\begin{lemma}
  If $U\in \Ar$ has two decompositions;
  \begin{center}
   $U=U_1U_{2}^{-1}=\overline{U}_1\overline{U}_{2}^{-1}$,
  \end{center}
  where $U_i$, $\overline{U}_i \in \pAr$ and $U_2$,
  $\overline{U}_2$ central in $\Ar$, then
  \begin{center}
    $U_1 <^+ U_2 \Longleftrightarrow \overline{U}_1 <^+
    \overline{U}_2$.
  \end{center}
\end{lemma}

\begin{proof}
  $U=U_1U_{2}^{-1}=\overline{U}_1\overline{U}_{2}^{-1}$ implies
  $U_1\overline{U}_2 \pe \overline{U}_1U_2$, ($\pe$ means equal in the monoid) since
  $U_2,\overline{U}_2$ central and $\pAr$ canonically injects in
  $\Ar$.\\
  If $U_1 <^+ U_2$ then
  \begin{equation*}
    \begin{split}
    &\Rightarrow U_1\overline{U}_2 <^+ U_2\overline{U}_2 \quad \text{since
      $<^+$ right-invariant}, \\
    &\Rightarrow U_1\overline{U}_2 <^+ \overline{U}_2U_2 \quad \text{since
      $U_2$ central}, \\
    &\Rightarrow \overline{U}_1U_2 <^+ \overline{U}_2U_2 \quad \text{since
      $U_1\overline{U}_2 \pe \overline{U}_1U_2$}, \\
    &\Rightarrow \overline{U}_1 <^+ \overline{U}_2,
    \end{split}
  \end{equation*}
    where the last implication follows from the fact that if
    $\overline{U}_{2} \le^+ \overline{U}_1$ then either: $(i)$
    $\overline{U}_2 = \overline{U}_1$, in which case $U=1$ and so
    $U_1=U_2$, a contradiction, or $(ii)$ $\overline{U}_2  <^+
    \overline{U}_1$, in which case $\overline{U}_2U_2  <^+
    \overline{U}_1U_2$. Again, a contradiction.

    The reverse implication follows by symmetry.
\end{proof}

This lemma shows that the following set is well defined:
\begin{equation*}
  \P = \{U \in \Ar : \text{ $U$ has decomposition } U=U_1U_{2}^{-1} \text{ where }
  U_2 <^+ U_1 \}.
\end{equation*}
It is an easy exercise to check that $\P$ is a positive cone in
$\Ar$ which contains $\P^+$: the positive cone in $\pAr$ with
respect to the order $<^+$.  Thus, the right-invariant order $<^+$
on $\pAr$ extends to a right-invariant order $<$ on $\Ar$, and we have shown the following.

\begin{theorem}
If $\Gam$ is a spherical-type Coxeter graph, then the Artin group $\Ar$ is right-orderable if and only if the corresponding Artin monoid $\pAr$ is right-orderable. \qed
\end{theorem}

\subsubsection{Reduction to type $E_8$}

Table \ref{tab:artin-injections} shows that every irreducible
spherical-type Artin group injects into one of type $A$, $D$, or $E$.
According to theorem \ref{artinro}, the Artin groups of type $A$ and $D$
are right-orderable. The Artin group of types
$E_6$ and $E_7$ naturally live inside $\A_{E_8}$, so it suffices
to show $\A_{E_8}$ is right-orderable, to conclude that all Artin groups are right-orderable.
 At this point in time it is
unknown whether $\A_{E_8}$ is right-orderable. As section
\ref{subsec:orderable-monoid} indicates it is enough to decide
whether the Artin monoid $\A_{E_8}^+$ is right-orderable.


 \vfill\eject {
 \baselineskip=14pt
 \bibliographystyle{alpha}
 \bibliography{artin}
}


\end{document}